\documentclass[a4paper,12pt]{article}
\usepackage{amssymb}
\usepackage{amsmath}
\usepackage{amsthm}
\usepackage{a4wide}
\usepackage{dsfont}            
\usepackage{mathrsfs}          
\usepackage{eufrak}             
\usepackage{accents}            
\usepackage{pstricks}
\usepackage{pst-plot}
\usepackage{pst-node}

\newcommand{\HH}{\ensuremath{\mathds H}}
\newcommand{\R}{\ensuremath{\mathds R}}						
\newcommand{\Rd}{\ensuremath{\mathds R^d}}				
\newcommand{\Nn}{\ensuremath{\mathds N}}		
\newcommand{\N}{\ensuremath{\mathds N_0}}				    
\newcommand{\Z}{\ensuremath{\mathds Z}}                      	
\newcommand{\Zd}{\ensuremath{\mathds Z^\mathrm{d}}}			
\newcommand{\wP}{\ensuremath{\mathds P}}				    
\newcommand{\E}{\ensuremath{\mathds E}}	                    
\newcommand{\one}{\ensuremath{\mathds 1}}                   
\newcommand{\domain}{\ensuremath{\mathcal{O}}}

\newcommand{\ijk}{\ensuremath{{i,j,k}}}                     
\newcommand{\jk}{\ensuremath{{j,k}}}
\newcommand{\Lp}{\ensuremath{L_p}}
\newcommand{\supp}{\ensuremath{\text{supp }}}

\newtheorem{thm}{Theorem}
\newtheorem{cor}[thm]{Corollary}

\theoremstyle{definition}
\newtheorem{ex}[thm]{Example}
\newtheorem{rem}[thm]{Remark}
\newtheorem{definition}[thm]{Definition}

\newcommand{\nnrm}[2]{\ensuremath{\| #1 \|_{#2}}}

\definecolor{felix}{rgb}{0.2,0.2,1.0} 
\definecolor{petru}{rgb}{0.7,0.1,0.1} 
\definecolor{dgruen}{rgb}{0.5,.7,.5}
\definecolor{ddblau}{rgb}{0.1,0,.55}
\definecolor{mylightGray}{gray}{0.95}

\newcommand{\fli}{\black}
\newcommand{\ilf}{\black}

\newcommand{\pci}{\black}
\newcommand{\icp}{\black}

\allowdisplaybreaks

\begin{document}

\title{Spatial Besov Regularity for\\ Stochastic Partial Differential Equations\\ on Lipschitz Domains}
\date{}
\maketitle
\begin{center}
\textsc{
Petru A. Cioica, Stephan Dahlke, Stefan Kinzel,
Felix Lindner,\\ Thorsten Raasch, Klaus Ritter, Ren{\'e} L. Schilling
}
\end{center}

\let\thefootnote\relax\footnotetext{\fli
This work has been supported by the Deutsche Forschungsgemeinschaft (DFG, grants DA 360/13-1, RI 599/4-1, SCHI 419/5-1) and a doctoral scholarship of the Philipps-Universit{\"a}t Marburg.
\ilf}

\begin{abstract}
\fli We use the scale of Besov spaces $B^\alpha_{\tau,\tau}(\domain),\;\alpha>0,\;1/\tau=\alpha/d+1/p$, $p$ fixed, to study the spatial regularity of the solutions of linear parabolic stochastic partial differential equations on bounded Lipschitz domains $\domain\subset\Rd$. The Besov smoothness determines the order of convergence that can be achieved by nonlinear approximation schemes. The proofs are based on a combination of weighted Sobolev estimates and characterizations of Besov spaces by wavelet expansions.\ilf
\end{abstract}
wavelet, weighted Sobolev space, nonlinear approximation, \fli adaptive numerical scheme \ilf

\section{Introduction}
In this paper, the spatial Besov regularity of the solutions of linear stochastic evolution equations on bounded Lipschitz domains is studied. We combine \fli regularity \ilf results by \textsc{Kim} \cite{Kim08} on stochastic partial differential equations (SPDEs, for short) on Lipschitz domains in terms of weighted Sobolev spaces with methods used in \textsc{Dahlke, DeVore} \cite{DD}, where the Besov regularity of (deterministic) elliptic equations on Lipschitz domains is investigated.  Our considerations are motivated by the question whether adaptive and other nonlinear approximation methods for the solutions of SPDEs on Lipschitz domains pay off in the sense that they  yield  better convergence rates than  uniform  methods. Thus referring to a numerical theme and combining concepts and methods from different areas and scientific communities, the article is addressed to readers of both worlds: stochastic analysis and numerical analysis.  Therefore, we give a \fli rather \ilf detailed account in the first part of the paper, emphasizing conceptual and notational clarity.

Our setting is as follows. On a finite interval $[0,T]\subset[0,\infty)$ let $(w^\kappa_t)_{t\in[0,T]},\;\kappa\in\mathds N=\{1,2,\ldots\}$, be independent, one-dimensional standard Brownian motions with respect to a filtration $(\mathcal F_t)_{t\in[0,T]}$ of $\sigma$-algebras on a complete probability space $(\Omega,\mathcal F,\wP)$. Throughout the paper we assume that $(\mathcal F_t)_{t\in[0,T]}$ is normal, i.e.\ the filtration satisfies the usual hypotheses, see, e.g.\ \cite[Section 3.3.]{DaPraZab}.
Let $\mathcal O\subset\Rd$ be a bounded Lipschitz domain.  We consider the model equation
 \begin{equation}\label{eq}
du=\sum_{\mu,\nu=1}^da^{\mu\nu}u_{x_\mu x_\nu}\,dt+\sum_{\kappa=1}^\infty g^\kappa\,dw^\kappa_t,\quad u(0,\,\cdot\,)=u_0,
\end{equation}
for $t\in[0,T]$ and $x\in \mathcal O$. Here $du$ is It{\^o}'s stochastic differential with respect to $t$, $(a^{\mu\nu})_{1\leq\mu,\nu\leq d}\in\R^{d\times d}$ is a strictly positive definite, symmetric matrix and the  coefficients  $g^\kappa,\;\kappa\in\mathds N,$ are random functions depending on $t$ and $x$ such that the mappings $\Omega\times[0,T]\ni(\omega,t)\mapsto g^\kappa(\omega,t,\,\cdot\,)$ are predictable processes with values in certain function spaces.  For details see  Section \ref{kimSection}.

Equation \eqref{eq} is understood in a weak  or distributional  sense, \frenchspacing{i.e.} $u$  is a solution of \eqref{eq}, if for all smooth and compactly supported test functions  $\varphi\in C_0^\infty(\domain)$ the equality
\begin{equation*}
\langle u(t,\,\cdot\,),\varphi\rangle=\langle  u_0  ,\varphi\rangle +\sum_{\mu,\nu=1}^d\int_0^t\langle a^{\mu\nu}u_{x_\mu x_\nu}(s,\,\cdot\,),\varphi\rangle\,ds+\sum_{\kappa=1}^\infty \int_0^t\langle g^\kappa(s,\,\cdot\,),\varphi\rangle\,dw^\kappa_s
\end{equation*}
holds for all $t\in[0,T]$ $\wP$-almost surely. Here and throughout the paper we write $\langle u,\varphi\rangle$ for the application of a distribution $u\in \mathcal D'(\domain)$ to a test function $\varphi\in C_0^\infty(\domain)$. \fli The \ilf existence and uniqueness of solutions of equation \eqref{eq}, respectively equation \eqref{eq''}  below,  within certain classes $\mathfrak H^\gamma_{p,\theta}(\domain,T)$ of stochastic processes has been shown in \cite{Kim08};  see also the earlier papers by  \textsc{Krylov}, \textsc{Lototsky} and \textsc{Kim}, \frenchspacing{e.g.} \cite{Kim04}, \cite{Kry99}, \cite{KryLot98}, \cite{Lot}. Roughly speaking, the classes $\mathfrak H^\gamma_{p,\theta}(\domain,T)$ are $L_p$-spaces of functions on $\Omega\times[0,T]$ with values in weighted Sobolev spaces $H^\gamma_{p,\theta-p}(\domain)$ that can be regarded as generalizations of the classical Sobolev spaces with  zero  Dirichlet boundary condition. Again we refer to Section \ref{kimSection} for precise  definitions.
Let us remark that in Examples \ref{ex1}, \ref{ex2} and \ref{ex3}, illustrating our Besov regularity result in Section \ref{resultSection}, the solution of equation \eqref{eq} in the class  $\mathfrak H^\gamma_{2,\theta}(\domain,T)$  coincides with the unique weak solution with \fli zero \ilf Dirichlet boundary condition in the sense of \textsc{Da Prato, Zabczyk} \cite{DaPraZab}, and hence can be represented by the well known stochastic  variation-of-constants formula
\begin{equation}\label{sgApproach}
u(t,\,\cdot\,)=e^{tA}u_0+\int_0^te^{(t-s)A}G(s)\,dW_s,\quad t\in[0,T].
\end{equation}
Here $(e^{tA})_{t\geq0}$ is the semigroup of contractions on $L_2(\domain)$  generated by the partial differential operator $A= \sum_{\mu,\nu=1}^d  a^{\mu\nu}\frac{\partial^2}{\partial x_\mu\partial x_\nu}$ with \fli zero \ilf Dirichlet boundary condition considered as an unbounded operator on $L_2(\domain)$, $(G(t))_{t\in[0,T]}$ is an operator-valued process and $(W_t)_{t\in[0,T]}$ is a cylindrical Wiener process on $\ell_2(\Nn)$, see Remarks \ref{remDaPra1} and \ref{remDaPra2} in Section \ref{kimSection}.

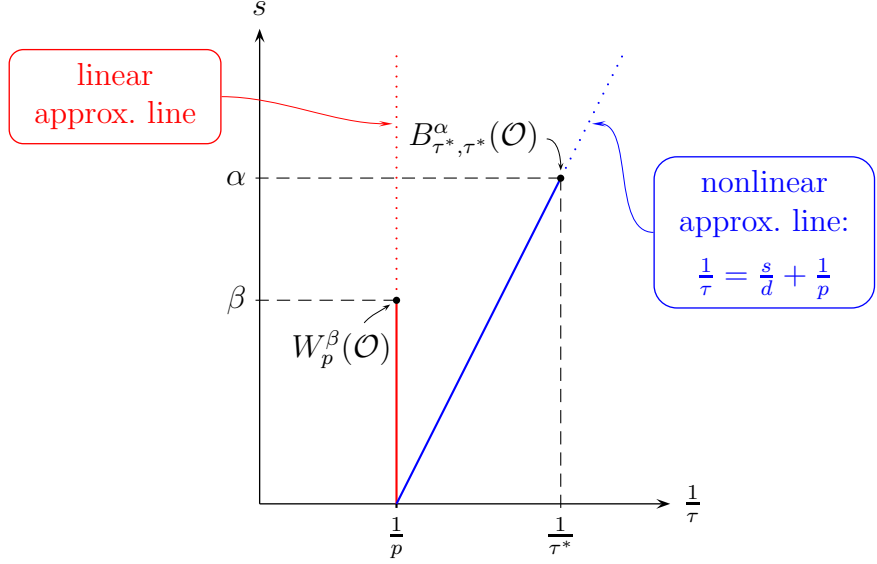
\begin{figure}

\psset{yunit=1.8,xunit=1.8}

\begin{center}
\begin{pspicture}(-1.5,0)(4,4)	

	\psaxes*[labels=none,ticks=none,linewidth=0.5pt,arrowscale=1.5]{->}(3,3.5)
	\uput[r](3,0){$\frac{1}{\tau}$}
	\uput[u](0,3.5){$s$}
	
	\psline(1,-0.03)(1,0)
	\psline[linecolor=red](1,0)(1,1.5)
	\psline[linecolor=red,linestyle=dotted](1,1.5)(1,3.3)
	
	 \uput{0.5cm}[l](0,3){\rnode[linecolor=red]{A}{\psframebox[fillstyle=solid,fillcolor=white,linewidth=.5pt,linecolor=red,framearc=.3]{\tabular{c} \red linear \\ \red approx. line\endtabular}}}
	\pnode(1,2.8){B}
	 \nccurve[nodesepB=.05,linecolor=red,linewidth=0.3pt,angleA=0,angleB=180]{->,arrowlength=2.5}{A}{B}

	\psline[linecolor=blue](1,0)(2.2,2.4)	
	\psline[linecolor=blue,linestyle=dotted](2.2,2.4)(2.65,3.3)
	 \uput{0.5cm}[r](2.6,2){\rnode[linecolor=blue]{G}{\psframebox[linewidth=.5pt,linecolor=blue, framearc=.3]{\tabular{c}{\blue  nonlinear} \\ {\blue approx. line:}\vspace{0.2cm}  \\ \blue $\frac{1}{\tau}=\frac{s}{d}+\frac{1}{p}$\endtabular}}}
	\pnode(2.4,2.8){H}
	 \nccurve[nodesepB=.05,linecolor=blue,linewidth=0.3pt,angleA=180,angleB=0]{->,arrowlength=2.5}{G}{H}

	\psdots[dotscale=.75](1,1.5)
	\uput[d](1,0){$\frac{1}{p}$}
	\uput[l](0,1.5){$\beta$}
	\psline[linewidth=0.4pt,linestyle=dashed](-0.03,1.5)(1,1.5)
	\uput{.4}[d](0.6,1.5){\rnode{C}{$W^\beta_p(\domain)$}}
	\pnode(1,1.5){D}
	\ncarc[linewidth=.3pt,nodesep=.1]{->}{C}{D}

	\psdots[dotscale=.75](2.2,2.4)
	\uput[l](0,2.4){$\alpha$}
	\psline[linewidth=.4pt,linestyle=dashed](-0.03,2.4)(2.2,2.4)
	\uput[d](2.2,0){$\frac{1}{\tau^*}$}
	\psline[linewidth=.4pt,linestyle=dashed](2.2,-0.03)(2.2,2.4)
	\uput{.4}[ul](2.2,2.4){\rnode{E}{$B^{\alpha}_{\tau^*,\tau^*}(\domain)$}}
	\pnode(2.2,2.4){F}
	\nccurve[linewidth=.3pt,angleA=0,angleB=90,nodesep=.1]{->}{E}{F}
\end{pspicture}

\end{center}
\caption[DeVoreTriebel001]{\begin{tabular}[t]{l}  Linear vs. nonlinear approximation \\  illustrated in \fli a \ilf\textsc{DeVore-Triebel} diagram.\end{tabular}}\label{fig:DeVoreTriebel001}

\end{figure}

As already mentioned, our motivation to study the Besov regularity of SPDEs is the theme of nonlinear approximation of the solution processes. For deterministic settings, a detailed  overview  of nonlinear approximation and an exposition of the characterization of its efficiency in terms of the Besov smoothness of the target functions can be found in \textsc{DeVore} \cite{DeV98}, see also \textsc{Cohen} \cite[\fli Chapters \ilf 3 and 4]{C03}.  Let us consider an example of approximation by wavelets in $L_p(\domain)$, the $L_p$-space of real-valued functions on $\domain,\;p\in(1,\infty)$. To this end, let  $\{\psi_\lambda\,:\,\lambda\in\nabla\}$  be a wavelet basis on $\domain$ and let $f\in L_p(\domain)$ be a target function which we want to approximate by functions $f_N\in L_p(\domain)$ belonging to certain approximation spaces  \fli$S_N$, where $N$ is the number of parameters used to describe the elements of $S_N$. \ilf We specify the index set of the wavelet basis by writing $\nabla=\bigcup_{j\geq j_0-1}\nabla_j$; the wavelets $\psi_\lambda,\;\lambda\in\nabla_j,\;j\geq j_0$, are those at scale levels $j\geq j_0$ respectively,  and $\psi_\lambda,\;\lambda\in\nabla_{j_0-1},$ are the scaling functions at the coarsest level $j_0\in\Z$.
In the case of  uniform  wavelet approximation up to a highest scale level $j_0-1+n,\;n\in\Nn,$  the approximation spaces are
\[S_N=S_{N(n)}=\Bigg\{\fli\sum_{j=j_0-1}^{j_0-1+n}\sum_{\lambda\in\nabla_j}c_\lambda\psi_\lambda\ilf\;:\;
c_\lambda\in\R,\;\lambda\in\nabla_j,\;j\in\{j_0-1,\ldots,j_0-1+n\}\Bigg\},\]
where $N=N(n)=|\bigcup_{j=j_0-1}^{j_0-1+n}\nabla_j|\in\Nn$ is the cardinality of the set of all indices up to scale level $j_0-1+n$. \fli Let $e_N(f)=\inf_{f_N\in S_N}\|f-f_N\|_{L_p(\domain)}$ be the corresponding approximation error measured in $L_p(\domain)$. It is well known that---under certain technical assumptions on the wavelet basis---the decay rate of $e_N(f)$ is linked to the $L_p$-Sobolev smoothness of the target function. \ilf More  precisely,  there exists an upper bound $r\in\Nn$ depending on the wavelet basis such that, for all $s\in[0,r]$,
\[f\in W^s_p(\domain)\quad\Longrightarrow\quad  e_N(f)\leq C\cdot N^{-s/d},\;N=N(n),\; n\in\Nn,\]
for some constant $C> 0$ which does not depend on $N$. The fractional order Sobolev spaces $W^s_p(\domain)$ are defined in the next section. One can also show \fli the converse \ilf
\[ \exists\,C>0\;\,\forall\,n\in\Nn\,:\,e_N(f)\leq C\cdot N^{-s/d},\;N=N(n)
\quad\Longrightarrow\quad f\in W^{s'}_p(\domain),\;s'<s.\]
If we consider instead best $N$-term approximation as a form of nonlinear approximation, the approximation spaces are
\[\Sigma_N=\Bigg\{\fli\sum_{\lambda\in\Lambda}c_\lambda\psi_\lambda\ilf\;:\;\Lambda\subset\nabla,\;|\Lambda|\leq N,\;c_\lambda\in\R,\;\lambda\in\Lambda\Bigg\},\]
$N\in\Nn$, and in this case the decay rate of the \fli error $\sigma_N(f):=\inf_{f_N\in\Sigma_N}\|f-f_N\|_{L_p(\domain)}$ \ilf is governed by the smoothness of $f$ measured in certain $L_\tau(\domain)$-norms, $\tau<p,$ which are weaker than the $L_p(\domain)$-norm: For all $\alpha\in[0,r]$,
\[ f\in B^\alpha_{\tau,\tau}(\domain),\quad\frac1\tau=\frac\alpha d+\frac1p\quad\Longrightarrow\quad\sigma_N(f)\leq C\cdot N^{-\alpha/d},\;N\in\Nn,\]
$B^\alpha_{\tau,\tau}(\domain)$ being a Besov space as defined in Section \ref{waveletsBesovSpaces}.
Therefore, if the target function $f$ belongs to $B^\alpha_{\tau^*,\tau^*}(\domain),\;1/\tau^*=\alpha/d+1/p$, for some $\alpha\in[0,r]$, and if in addition $\beta:=\sup\{s\in\R\,:\,f\in W^s_p(\domain)\}<\alpha$, then the convergence rate of  uniform wavelet approximations  is inferior to the convergence rate of  the best $N$-term wavelet  approximation. The latter can be considered as a benchmark for the convergence rate of adaptive numerical algorithms,  see \fli \cite{CDD1}, \cite{CDD2}, \cite{DDD}. This \ilf situation is illustrated in Figure \ref{fig:DeVoreTriebel001}, where each point $(1/\tau,s)$ represents the smoothness spaces of functions with ``$s$ derivatives in $L_{\tau}(\domain)$''.  Note that the nonlinear approximation line $\{(1/\tau,s)\in[0,\infty)^2\,:\,1/\tau=s/d+1/p\}$ is also the Sobolev embedding line\fli. For bounded domains, \ilf all spaces left to this line as well as the spaces $B^s_{\tau,\tau}(\domain)$ on the line are continuously embedded in $L_p(\domain)$. 


Let us return to equation \eqref{eq} and assume that the solution $u=u(\omega,t,x),\;(\omega,t,x)\in\Omega\times[0,T]\times\domain,$ vanishes  on  the boundary $\partial\domain$, satisfying a  zero Dirichlet boundary condition. It is clear that the smoothness of $x\mapsto u(\omega,t,x)$ depends on the smoothness of the mappings $x\mapsto g^\kappa(\omega,t,x),\;\kappa\in\Nn$. However, even if the spatial smoothness of the $g^\kappa$ is high, the Sobolev smoothness of $x\mapsto u(\omega,t,x)$ can be additionally limited by singularities of the spatial derivatives of $u$ at  the boundary  of $\domain$, due to the  zero  Dirichlet boundary condition \fli and the shape of the domain\ilf.
Such corner singularities are \fli typical  examples for the fact that \ilf the spatial $L_p$-Sobolev regularity of $u$ \fli may be \ilf exceeded by the regularity in the scale of Besov spaces $B^\alpha_{\tau,\tau}(\domain),\fli\;1/\tau=\alpha/d+1/p\ilf$. In this paper, we present a result on the spatial Besov regularity of the solution $u$ to equation \eqref{eq} which has the following structure: If
\[u\in L_p(\Omega\times[0,T],\mathcal P,\wP\otimes\lambda;\,W^s_p(\domain))\]
and if the  functions  $g^\kappa,\;\kappa\in\Nn$, are sufficiently regular, then
\[u\in L_\tau(\Omega\times[0,T],\mathcal P,\wP\otimes\lambda;\,B^\alpha_{\tau,\tau}(\domain))\] for certain $\alpha>s$ and \fli $1/\tau=\alpha/d+1/p$\ilf. Here $\mathcal P$ is the predictable $\sigma$-algebra w.r.t.\  the filtration $(\mathcal F_t)_{t\in[0,T]}$ and $\lambda$ denotes Lebesgue measure on $[0,T]$.  This result is important  for the theoretical foundation of adaptive numerical methods for the approximation of $u$.
The proof is based on a wavelet expansion of an extension of $\domain\ni x\mapsto u(\omega,t,x)$ to $\Rd$, which allows \fli us to estimate \ilf the $B^\alpha_{\tau,\tau}(\domain)$-norm in terms of the wavelet coefficients. We apply a strategy similar to the one used in \textsc{Dahlke, DeVore} \cite{DD}, where the Besov regularity of (deterministic) elliptic equations on Lipschitz domains is investigated with the help of an estimate of weighted Sobolev norms of harmonic functions. Our substitute for the latter is an estimate of weighted Sobolev norms of the solution of equation \eqref{eq} provided by \textsc{Kim} \cite{Kim08}. 

There \fli exists an extensive literature on \ilf the Besov regularity of SPDEs. In general, however, the assumptions on the domain and the scale of parameters considered do not fit into our setting. To mention an example, the semigroup approach to SPDEs of \textsc{Da Prato, Zabczyk} \cite{DaPraZab}, which is placed in a Hilbert space framework, has been generalized to M-type $2$ Banach spaces by \textsc{Brze{\'z}niak} \cite{Brz95}, \cite{Brz97}, for the purpose of gaining better H{\"o}lder regularity results. Roughly speaking, the operator $A$ appearing in equation $\eqref{sgApproach}$ is considered as the generator of a semigroup on $L_p(\domain)$ for some $p\geq2$, and the stochastic integral in \eqref{sgApproach} is considered as an stochastic integral in an interpolation space $X$ between $L_p(\domain)$ and $D(A)\subset L_p(\domain)$, the domain of $A$, realizing a zero Dirichlet boundary condition. If $\partial\domain$ is sufficiently smooth, then $D(A)=W^2_p(\domain)\cap\accentset{\circ}{W}^1_p(\domain)$ and $X\subseteq B^s_{p,2}(\domain)$ for some $s\in[0,2]$. In this situation, \fli the Sobolev embedding theorem leads to H{\"o}lder regularity results, and these results become better for large $p$\ilf. With the help of a theory of stochastic integration in wider classes of Banach spaces, this approach has been generalized in the works of \textsc{Van Neerven, Weis, Dettweiler} and \textsc{Veraar}, see, e.g.\ \fli \cite{DVNW}, \cite{VNW05}, \cite{VNVW07}, \cite{VNVW08}\ilf, compare also \textsc{Brze{\'z}niak, Van Neerven} \cite{BrzVN00}. In contrast to these works the problem \fli considered here \ilf is of a different nature. Firstly, we are explicitly interested in domains with  non-smooth boundary. For polygonal non-convex domains, it is well known that $W^2_2(\domain)\cap\accentset{\circ}{W}^1_2(\domain)\varsubsetneq D(A)$, where $D(A):=\{u\in\accentset{\circ}{W}^1_2(\domain):Au\in L_2(\domain)\},\; A=\Delta=\sum_{\mu=1}^d\frac{\partial^2}{\partial x_\mu^2}$, \fli see \textsc{Grisvard} \cite{Gri85}, \cite{Gri92}, and for more general Lipschitz domains see \textsc{Jerison, Kenig} \cite{JK}\ilf. Secondly, we are interested in the special scale $B^\alpha_{\tau,\tau}(\domain),\;1/\tau=\alpha/d+1/p,\;\fli\tau>0\ilf,\;p$ fixed, including in particular spaces which are no Banach spaces but quasi-Banach spaces. The parameter $\tau$ decreases if $\alpha$ increases and $B^\alpha_{\tau,\tau}(\domain)$ fails to be a Banach space for $\tau<1$.
\fli While our methods work in this setting, any direct approach requires (at least!) a fully-fledged theory of stochastic integration in quasi-Banach spaces which is not yet available.\ilf

\fli Let \ilf us emphasize that our result can be extended to more general linear equations of the type
\begin{equation}\label{eq''}
\left.
\begin{aligned}
&du=\sum_{\mu,\nu =1}^{d}\left(a^{\mu\nu}u_{x_\mu x_\nu}+b^\mu u_{x_\mu}+cu+f\right)dt +
\sum_{\kappa=1}^{\infty}\left(\sum_{\mu=1}^d\sigma^{\mu\kappa}u_{x_\mu}+\eta^\kappa u+g^\kappa\right)dw^\kappa_t, \\
&u(0,\,\cdot\,)=u_0,
\end{aligned}\;
\right\}
\end{equation}
including, in particular, the case of multiplicative noise. Here the coefficients  $a^{\mu\nu}$, $b^\mu$, $c$, $\sigma^{\mu\kappa}$, $\eta^\kappa$  and the free terms $f$ and $g^\kappa$ are random functions depending on $t$ and $x$. This extension is possible because one of our main  tools,  the weighted Sobolev norm estimate of Corollary \ref{cor} in Section \ref{kimSection}, holds for equations of type \eqref{eq} as well as for equations of type \eqref{eq''}. \fli Since \ilf this mainly adds notational complications, we will focus on equation \eqref{eq} and refer to Appendix \ref{AppB} for a short account of how to treat equations of type \eqref{eq''}.
\bigskip

The paper is organized as follows: In Section 2 we  collect  the notations, definitions and preliminary results needed later on. Some general notations are introduced in Section \ref{notations}. Section \ref{waveletsBesovSpaces} provides the necessary facts  on  Besov spaces and wavelet decompositions. In Section \ref{kimSection} a short introduction to the general $L_p$-theory of SPDEs on Lipschitz domains due to \textsc{Kim} \cite{Kim08} is given, including definitions of the already mentioned spaces $H^\gamma_{p,\theta-p}(\domain),\;\mathfrak H^\gamma_{p,\theta}(\domain,T)$. Finally, in Section \ref{resultSection} the Besov regularity result (Theorem \ref{result}) is stated and proved, and some concrete examples for an application of the result are given.

\section{Preliminaries}

\subsection{Some notations and conventions}\label{notations}
In this and the next subsection $\domain\subseteq\Rd$ can be an arbitrary (not necessarily bounded) Lipschitz domain. A domain is called Lipschitz if  each point on the boundary $\partial\domain$ has a \fli neighbourhood \ilf whose intersection with the boundary---after relabeling and reorienting the coordinate axes if necessary---is the graph of a Lipschitz function

By $\mathcal D'(\domain)$ we denote the space of Schwartz distributions on $\domain$. \fli  If not explicitly stated otherwise, all function spaces or spaces of distributions are meant to be spaces of real-valued functions or distributions. \ilf If $f\in\mathcal D'(\domain)$ is a generalized function and $\alpha=(\alpha_1,\ldots,\alpha_d)\in\N^d$ is a multi-index, we write $D^\alpha f=\frac{\partial^{|\alpha|}f}{\partial x_1^{\alpha_1}\ldots\partial x_d^{\alpha_d}}$ for the corresponding derivative w.r.t. $x=(x_1,\ldots,x_d)\in\domain$, where $|\alpha|=\alpha_1+\ldots+\alpha_d$. As in equations \eqref{eq} and \eqref{eq''} we also use the notation $f_{x_\mu x_\nu}=\frac{\partial^{2}f}{\partial x_\mu\partial x_\nu},\;f_{x_\mu}=\frac{\partial f}{\partial x_\mu}.$ For $m\in\N$, $D^mf=\{D^\alpha f\,:\, |\alpha|=m\}$ is the set of all $m$-th order derivatives of $f$
which is identified with an  $\R^{\tbinom {d+m-1}m}$-valued  distribution. Given $p\in[1,\infty)$ and $m\in\N$, $W^m_p(\domain)$ denotes the classical Sobolev space consisting of all (equivalence classes of) measurable functions $f:\domain\to\R$ such that $\|f\|_{W^m_p(\domain)}=\|f\|_{L_p(\domain)}+|f|_{W^m_p(\domain)}=(\int_\domain |f(x)|^p\,dx)^{1/p}+\sum_{|\alpha|=m}(\int_\domain|D^\alpha f(x)|^p\,dx)^{1/p}$ is finite. For $p\in(1,\infty)$ and $s\in(m,m+1),\;m\in\N$, we define the fractional order Sobolev space $W^s_p(\domain)$ to be the Besov space $B^s_{p,p}(\domain)$ introduced in the next subsection.  (This scale of fractional order Sobolev spaces can also be obtained by real interpolation of $W^n_p(\domain),\;n\in\N$. One can show that $W^n_2(\domain)=B^n_{2,2}(\domain)$ for all $n\in\Nn$ and $W^n_p(\domain)\subset B^n_{p,p}(\domain)$ for all $n\in\Nn,\;p>2$, \frenchspacing{see, e.g.} \textsc{Triebel} \cite[Remark 2.3.3/4 and Theorem 4.6.1.(b)]{Tri78} together with \textsc{Dispa} \cite{Dis}.)
Given any countable index set $\mathcal J$, the space of $p$-summable sequences indexed by $\mathcal J$ is denoted by $\ell_p=\ell_p(\mathcal J)$ and $|\,\cdot\,|_{\ell_p}$ is the respective norm. Usually we have $\ell_p=\fli\ell_p(\Nn)\ilf$ but, for instance we may also use the notation $|D^m f(x)|_{\ell_p}^p=\sum_{|\alpha|=m}|D^\alpha f(x)|^p$ for $f\in W^m_p(\domain)$.

Given a distribution $f\in D'(\domain)$ and a smooth and compactly supported test function $\varphi\in C_0^\infty(\domain)$, we write  $\langle f,\varphi\rangle$  for the application of $f$ to $\varphi$. If $H$ is a Hilbert space, then $\langle\,\cdot\,,\,\cdot\,\rangle_H$ denotes the inner product in $H$. Given another Hilbert space $U$, we denote by $L_{(\text{HS})}(H,U)$ and $L_{(\text{nuc})}(H,U)$ the spaces of Hilbert-Schmidt operators and nuclear operators from $H$ to $U$ respectively, see, e.g. \fli \textsc{Pietsch} \cite[Sections 6 and 15]{Pie} \ilf or \textsc{ Da~Prato, Zabczyk} \cite[Appendix C]{DaPraZab} for definitions. We also abbreviate $L_{(\text{HS})}(H)=L_{(\text{HS})}(H,H)$ and $L_{(\text{nuc})}(H)=L_{(\text{nuc})}(H,H)$. $\mathcal M_T^{2,c}(H,(\mathcal F_t))$ is the space of continuous, square integrable, $H$-valued martingales with respect to the filtration $(\mathcal F_t)_{t\in[0,T]}$. For $\Omega\times[0,T]$ we use the  shorthand notation $\Omega_T$ and
$$\mathcal P=\sigma\big(\{]s,t]\times F_s\,:\,0\leq s<t\leq T,\,F_s\in\mathcal F_s\}\cup\{\{0\}\times F_0\,:\,F_0\in\mathcal F_0\}\big)$$
is the predictable $\sigma$-algebra. \fli $\wP\otimes\lambda$ is the product measure of the probability measure $\wP$ on $(\Omega,\mathcal F)$ and Lebesgue measure $\lambda$ on $([0,T],\mathcal B([0,T]))$, where $\mathcal B([0,T])$ denotes the Borel $\sigma$-algebra on $[0,T]$. Given any measure space $(A,\mathcal A,m)$, any (quasi-)normed space $B$ with (quasi-)norm $\|\cdot\|_B$ and any summability index $p>0$, we denote by $L_p(A,\mathcal A,m;B)$ the $L_p$-space of all strongly measurable functions $u:A\to B$ whose (quasi-)norm $\|u\|_{L_p(A,\mathcal A,m;B)}:=\left(\int_A\|u(z)\|_B^p\,m(dz)\right)^{1/p}$ is finite.\ilf

All equalities of random variables or random (generalized) functions appearing in this paper are meant to be $\wP$-almost sure equalities.
Throughout the paper, $C$ denotes a positive constant which may change its value  from line to line.

\subsection{Besov spaces and wavelet decompositions}\label{waveletsBesovSpaces}

In this section we give  the  definition of Besov spaces and describe their characterization in terms of wavelets. Our standard reference in this context is the monograph of \textsc{Cohen} \cite{C03}.

For a function $f:\domain \to \R$ and a natural number $n\in\fli\Nn\ilf$ let
\begin{equation*}
\Delta_h^n f (x) :=
\prod_{i=0}^n\one_\domain( x+ih ) \cdot \sum\limits_{j=0}^n {n \choose j} (-1)^{n-j}\, f(x+jh)
\end{equation*}
be the $n$-th difference of $f$ with step $h\in\Rd$.
For $p\in\left(0,\infty\right)$ the modulus of smoothness is given by
\[
\omega^n (t,f)_p := \sup_{|h|<t} \, \| \,
\Delta_h^n f \,  \|_{L_p (\domain)}\, , \qquad t>0 \, .
\]
One approach to introduce Besov spaces is the following.

\begin{definition}\label{BesovSpaces}
 Let $s, p, q \in\left(0,\infty\right)$ and $n\in\Nn$ with $n>s$.
Then $B^s_{p,q}(\domain)$ is the collection of all functions  $f \in L_p (\domain) $ such that
\[
| \, f \, |_{B^s_{p,q}(\domain)} := \bigg( \int_0^\infty \Big[
t^{-s} \, \omega^n (t,f)_p\Big]^q \frac{dt}{t}\bigg)^{1/q} <\infty.
\]
These classes are equipped with a  (quasi-)norm  by taking
\[
\| \, f \, \|_{B^s_{p,q}(\domain)}
:= \| \, f\, \|_{L_p(\domain)} +  | \, f \, |_{B^s_{p,q}(\domain)}\, .
\]
\end{definition}

\begin{rem}
For a more general definition of Besov spaces, including the cases where $p,q=\infty$ and $s<0$ see, \frenchspacing{e.g.} \textsc{Triebel} \cite{TriIII}.
\end{rem}


We want to describe $B^s_{p,q}(\Rd)$ by  means of wavelet expansions.  To this end let $\varphi$ be a scaling function of tensor product type on $\Rd$ and let $\psi_i$, $i=1, \ldots, 2^d-1$, be corresponding multivariate  mother  wavelets, such that, for a given $r\in\Nn$ and some $N>0$, the following locality, smoothness and vanishing moment conditions hold. For all $i=1, \ldots, 2^d-1$,
 \begin{align}
&\supp \varphi,\,\supp\psi_i\subset [-N,N]^d,\label{wl1}\\
&\varphi,\,\psi_i \in C^r(\Rd),\label{wl2}\\
&\int x^\alpha \, \psi_i (x)\, dx=0 \quad\text{ for all $\alpha\in\N^d$ with $|\alpha|\le r$}.\label{wl3}
\end{align}
We \fli assume \ilf that
\begin{align*}
 \big\{\varphi_k,\psi_{i,j,k}\,:\, (\ijk)\in\{1,\cdots,2^d-1\}\times\N\times\Z^d\big\}
\end{align*}
is a Riesz basis of $L_2(\Rd)$,  where we use  the standard abbreviations for dyadic shifts and dilations of the scaling function and the corresponding wavelets
 \begin{align}
\varphi_k(x)	&:=\varphi(x - k),\;x\in\Rd ,&&\text{for $k\in\Z^d$, and} \label{wl4}\\
\psi_{i,j,k}(x)	&:=2^{jd/2}\psi_i(2^jx-k),\;x\in\Rd, &&\text{for $(\ijk)\in\{1,\cdots,2^d-1\}\times\N\times\Z^d$.}\label{wl5}
\end{align}
Further, we assume that \fli there exists \ilf a dual Riesz basis satisfying the same requirements. More precisely, there exist functions
$\widetilde{\varphi} $ and
$\widetilde{\psi}_i$, $ i=1, \ldots , 2^d-1$,
 such that conditions \eqref{wl1}, \eqref{wl2} and \eqref{wl3} hold if $\varphi$ and $\psi$ are replaced by $\widetilde{\varphi} $ and
$\widetilde{\psi}_i$, and such that the biorthogonality relations
\[
\langle \widetilde{\varphi}_k, \psi_{i,j,k} \rangle = \langle \widetilde{\psi}_{i,j,k},
 \varphi_k \rangle  = 0\, ,	\quad
\langle \widetilde{\varphi}_k, \varphi_{\ell} \rangle  = \delta_{k,\ell}, \quad
\langle \widetilde{\psi}_{i,j,k}, \psi_{u,v,\ell} \rangle  = \delta_{i,u}\, \delta_{j,v}\, \delta_{k,\ell}\, ,
\]
are fulfilled.
Here we use  analoguous  abbreviations to \eqref{wl4} and \eqref{wl5} for the dyadic shifts and dilations of $\widetilde{\varphi} $ and
$\widetilde{\psi}_i$ ,  and $\delta_{k,l}$ denotes the Kronecker symbol. We refer to \textsc{Cohen} \cite[Chapter 2]{C03} for the construction of biorthogonal wavelet bases, see also \textsc{Daubechies} \cite{Daub} and \textsc{Cohen, Daubechies, Feauveau} \cite{CDF}. To  keep notation simple,  we will write
\begin{equation*}
 \psi_{i,j,k,p} := 2^{jd(1/p-1/2)}\psi_\ijk \qquad \text{ and } \qquad \widetilde{\psi}_{i,j,k,p'}:= 2^{jd(1/{p'}-1/2)}\widetilde{\psi}_\ijk,
\end{equation*}
for the $L_p$-normalized wavelets and the correspondingly modified duals, with $p':=p/(p-1)$ if $p\in(0,\infty),\;p\neq 1,$ and  $p':=\infty,\;1/p':=0$ if $p=1$.

The following theorem shows how Besov spaces can be described by decay properties of the wavelet coefficients, \fli if \ilf the parameters fulfil certain conditions.

\begin{thm}\label{BesovChar01}
Let $p,q\in\left(0,\infty\right)$ and $s>\max\left\{0,d\left(1/p-1\right)\right\}$. Choose $r\in\Nn$ such that $r>s$ and construct a  biorthogonal  wavelet Riesz basis as described above.
Then a locally integrable function $f:\Rd\to\R$ is in the Besov space $B^s_{p,q}(\Rd)$ \fli if, and only if,\ilf
\begin{equation}\label{BesovZerlegung}
f=\sum_{k\in\Zd}\langle f,\widetilde{\varphi}_k\rangle\,\varphi_k
	+ \sum_{i=1}^{2^{d}-1}\sum_{j\in\N}\sum_{k\in\Zd}\langle f,\widetilde{\psi}_{\ijk,p'}\rangle\,\psi_{\ijk,p}
\end{equation}
(convergence in $\mathcal D'(\Rd)$) with
\begin{equation}\label{BesovNormDiscrete}
\Big(\sum_{k\in\Zd}|\langle f,\widetilde{\varphi}_k\rangle|^p\Big)^{1/p}
	+ \Big(
	\sum_{i=1}^{2^{d}-1}\sum_{j\in\N}2^{j s q}
	\Big(\sum_{k\in\Zd}|\langle f,\widetilde{\psi}_{\ijk,p'}\rangle|^p\Big)^{q/p}\Big)^{1/q}	
	<	\infty,
\end{equation}
and \eqref{BesovNormDiscrete}  is an equivalent (quasi-)norm  for $B^s_{p,q}(\Rd)$.
\end{thm}

\begin{rem}
A proof of this theorem for the case $p\geq1$ can be found in \textsc{Meyer} \cite[\textsection10 of Chapter 6]{Mey92}. For the general case see for example \textsc{Kyriazis} \cite{Kyr96} or \textsc{Cohen} \cite[Theorem 3.7.7]{C03}. Of course, if \eqref{BesovNormDiscrete} holds then the infinite sum in \eqref{BesovZerlegung} converges also in $B^s_{p,q}(\Rd)$.  If $s>\max\left\{0,d\left(1/p-1\right)\right\}$ we have  the embedding $B_{p,q}^s(\Rd)\subset L_u(\Rd)$ for some $u>1$, see, \frenchspacing{e.g.} \textsc{Cohen} \cite[Corollary 3.7.1]{C03}.
\end{rem}

Let us now fix a value $p\in\left(1,\infty\right)$ and consider the scale of Besov spaces $B^s_{\tau,\tau}(\Rd)$, $\fli1/\tau=s/d+1/p\ilf$, $s>0$.
A simple computation gives the following result.

\begin{cor}\label{BesovChar02}
Let $p\in\left(1,\infty\right)$, $s>0$ and $\tau\in\R$ such that  $1/\tau=s/d+1/p$. Choose $r\in\Nn$ such that $r>s$ and construct a  biorthogonal  wavelet Riesz basis as described above.
Then a locally integrable function $f:\Rd\to\R$ is in the Besov space $B^s_{\tau,\tau}(\Rd)$ \fli if, and only if,\ilf
\begin{equation}\label{BesovZerlegung02}
f=\sum_{k\in\Zd}\langle f,\widetilde{\varphi}_k\rangle\,\varphi_k
	+ \sum_{i=1}^{2^{d}-1}\sum_{j\in\N}\sum_{k\in\Zd}\langle f,\widetilde{\psi}_{\ijk,p'}\rangle\,\psi_{\ijk,p}
\end{equation}
(convergence in $\mathcal D'(\Rd)$) with
\begin{equation}\label{BesovNormDiscrete02}
\Big(\sum_{k\in\Zd}|\langle f,\widetilde{\varphi}_k\rangle|^\tau\Big)^{1/\tau}
	+ \Big(
	\sum_{i=1}^{2^{d}-1}\sum_{j\in\N}
	\sum_{k\in\Zd}|\langle f,\widetilde{\psi}_{\ijk,p'}\rangle|^\tau\Big)^{1/\tau}	
	<	\infty	\, ,
\end{equation}
and \eqref{BesovNormDiscrete02}  is an equivalent (quasi-)norm  for $B^s_{\tau,\tau}(\Rd)$.
\end{cor}

%

\subsection{SPDEs on Lipschitz domains and weighted Sobolev spaces}\label{kimSection}

\fli From now on, let $\domain\subset\Rd$ be a bounded Lipschitz domain and $d\geq2$.\ilf

We have already mentioned corner singularities as  typical examples  where the \fli regularity of a function on $\domain\subset\Rd$ in the Besov scale \ilf $B^\alpha_{\tau,\tau}(\domain),\;1/\tau=\alpha/d+1/p,\;\alpha>0$,  can  exceed the regularity in the Sobolev scale $W^s_p(\domain),\;s>0$. This reflects the sparsity of the large wavelet coefficients of \fli such a \ilf function (given a wavelet basis on the domain $\domain$). A general way to deal with smoothness regardless of certain singularities at the boundary is  to use  weighted Sobolev spaces, where the weight function is a power of the distance to the boundary. The $L_p$-theory of SPDEs on Lipschitz domains by \textsc{Kim} \cite{Kim08} is based on spaces of this type, namely
the weighted Sobolev spaces $H^\gamma_{p,\theta}(\mathcal O),\;p\in(1,\infty),\;\theta,\,\gamma\in\mathds R,$ introduced in \textsc{Lototsky} \cite{Lot}. They are defined in terms of the Bessel-potential spaces
 \[H^\gamma_p(\R^d)=\{u\in \mathcal S'(\R^d)\,:\,\|u\|_{H^\gamma_p(\R^d)}=\|(1-\Delta)^{\gamma/2}u\|_{L_p(\R^d)}<\infty\}.\]
Here, $\mathcal S'(\Rd)\subset\mathcal D'(\Rd)$ is the space of (real valued) tempered distributions and $(1-\Delta)^{\gamma/2}:\mathcal S'(\Rd)\to\mathcal S'(\Rd)$ is the pseudo-differential operator with symbol $\Rd\ni\xi\mapsto(1+|\xi|^2)^{\gamma/2}$, \frenchspacing{i.e.} $(1-\Delta)^{\gamma/2}u=\mathcal F^{-1}\big((1+|\xi|^2)^{\gamma/2}\mathcal Fu\big)$,
where  $\mathcal F$ denotes the Fourier transform on the (complex valued) tempered distributions.

For $x\in \mathcal O$ we write $\rho(x):=\text{dist}(x,\partial \mathcal O)$ for the distance between $x$ and the boundary of the domain $\mathcal O$. Fix $c>1,\;k_0>0$ and for $n\in\Z$ consider the subsets $\mathcal O_n$ of $\mathcal O$ given by
\[\mathcal O_n:=\{x\in\mathcal O\,:\,c^{-n-k_0}<\rho(x)<c^{-n+k_0}\}.\]
Let $\zeta_n,\;n\in\Z,$ be non-negative functions satisfying $\zeta_n\in C_0^\infty(\mathcal O_n),\;\sum_{n\in\Z}\zeta_n(x)=1$ and  $|D^m\zeta_n(x)|\leq C\cdot c^{mn}$ for all $n\in\Z,\;m\in\N,\;x\in\mathcal O$, and a constant $C>0$ that does not depend on $n,\;m$ and $x$. The functions $\zeta_n$ can be constructed by mollifying the indicator functions of the sets $\domain_n$, see, e.g.\ \textsc{H{\"o}rmander} \cite[Section 1.4]{Hoer90}. \fli If $\domain_n$ is empty we set $\zeta_n\equiv0$. For $u\in\mathcal D'(\domain)$ $\zeta_nu$ is a distribution on $\domain$ with compact support which can be extended by zero to $\R^d$\ilf. This extension is a tempered distribution, \frenchspacing{i.e.} $\zeta_nu\in\mathcal S'(\Rd)$.
\begin{definition}
Let $\zeta_n,\;n\in\Z,$ be as above and $p\in(1,\infty),\;\theta,\,\gamma\in\mathds R$. Then
\[H^\gamma_{p,\theta}(\mathcal O):=\Big\{ u\in\mathcal D'(\mathcal O)\,:\,\|u\|_{H^\gamma_{p,\theta}(\mathcal O)}^p:=\sum_{n\in\Z}c^{n\theta}\|\zeta_{-n}(c^n\,\cdot\,)u(c^n\,\cdot\,)\|_{H^\gamma_p(\Rd)}^p
<\infty\Big\}.\]
\end{definition}
According to \textsc{Lototsky} \cite{Lot} this definition is independent of the specific choice of $c,\;k_0$ and $\zeta_n,\;n\in\N,$ in the sense that one gets equivalent norms.
If $\gamma=m\in\N$ then the spaces can be characterized as
\begin{align*}
H^0_{p,\theta}(\mathcal O)&=L_{p,\theta}(\mathcal O):= \Lp(\mathcal O,\rho(x)^{\theta-d}dx),\\
H^m_{p,\theta}(\mathcal O)&=\left\{u\,:\,\rho^{|\alpha|}D^\alpha u\in L_{p,\theta}(\mathcal O)\text{ for all }\alpha\in\N^d\text{ with } |\alpha|\leq m\right\},
\end{align*}
and one has the norm equivalence
\begin{equation}\label{ne}
C^{-1} \|u\|^p_{H_{p,\theta}^m(\mathcal O)}\;\leq\;\sum_{\alpha\in\N^d,\,|\alpha|\leq m}\int_\mathcal O\left|\rho(x)^{|\alpha|}
D^\alpha u(x)\right|^p\rho(x)^{\theta-d}\,dx\;\leq\;C\|u\|^p_{H_{p,\theta}^m(\mathcal O)}.
\end{equation}

Analogous notations are used for $\ell_2=\ell_2(\mathds N)$-valued functions $g=(g^\kappa)_{\kappa\in\mathds N}$. For $p\in(1,\infty),\;\theta,\,\gamma\in\R$ and $\zeta_n,\;n\in\Z$, as above
\begin{align*}
H^\gamma_p(\Rd;\ell_2)&:=\Big\{g\in(\mathcal S'(\Rd))^\Nn\,:\,(1-\Delta)^{\gamma/2}g^\kappa\in L_p(\Rd) \text{ for all } k\in\Nn \text{ and }\\
&\qquad\qquad\qquad\qquad\quad\;\,\|g\|_{H^\gamma_p(\Rd;\ell_2)}:=
\big\|\big|\big((1-\Delta)^{\gamma/2}g^\kappa\big)_{\kappa\in\Nn}\big|_{\ell_2}\big\|_{L_p(\Rd)}<\infty\Big\},\\
H^\gamma_{p,\theta}(\domain;\ell_2)&:=\Big\{ g\in(\mathcal D'(\domain))^\Nn\,:\,
\|g\|_{H^\gamma_{p,\theta}(\domain;\ell_2)}^p:=\sum_{n\in\Z}c^{n\theta}
\|\zeta_{-n}(c^n\,\cdot\,)g(c^n\,\cdot\,)\|^p_{H^\gamma_p(\Rd;\ell_2)}<\infty\Big\}.
\end{align*}
\begin{rem}\label{Kufner}
{\bf(a)} One can consider the spaces $H^\gamma_{p,\theta}(\domain)$ as generalizations of the classical Sobolev spaces on $\domain$ with zero boundary conditions. For $\gamma=m\in\N$ we have the identity
\[H^m_{p,d-mp}(\domain)=\accentset{\circ}{W}^m_p(\domain),\]
and the norms in both spaces are equivalent, see Theorem \frenchspacing{9.7.} in \textsc{Kufner} \cite{Kuf}. Here $\accentset{\circ}{W}^m_p(\domain)$ is the closure of $C_0^\infty(\domain)$ in the classical Sobolev space $W^m_p(\domain)$.

{\bf(b)} Note that, in contrast to the spaces $W^s_p(\domain)=B^s_{p,p}(\domain),\;s\in(m,m+1),\;m\in\N,$ which can be regarded as \emph{real} interpolation spaces of the classical Sobolev spaces $W^m_p(\domain),\;m\in\N$ (see, \frenchspacing{e.g.} \textsc{Triebel} \cite[\fli Section \ilf 1.11.8]{TriIII} \fli and \textsc{Dispa}\cite{Dis}\ilf), the spaces $H^\gamma_{p,\theta}(\domain),\;\gamma\in(m,m+1),\;m\in\N$, are \emph{complex} interpolants of the respective integer smoothness spaces (\frenchspacing{cf.} \textsc{Lototsky} \cite[\fli Proposition \ilf 2.4]{Lot}).
\end{rem}

 We can now  define spaces of stochastic processes and random functions in terms of the weighted Sobolev spaces introduced above. 
\begin{definition}
For $\gamma,\;\theta\in\R$ and $p\in(1,\infty)$ we set
\begin{align*}
\HH^\gamma_{p,\theta}(\mathcal O,T)&:=L_p\left(\Omega_T,\mathcal P,\wP\otimes\lambda;\,H^\gamma_{p,\theta}(\mathcal O)\right),\\
\HH^\gamma_{p,\theta}(\mathcal O,T;\ell_2)&:=L_p\left(\Omega_T,\mathcal P,\wP\otimes\lambda;\,H^\gamma_{p,\theta}(\mathcal O;\ell_2)\right),\qquad\qquad\qquad\qquad\\
U^\gamma_{p,\theta}(\mathcal O)&:=L_p\big(\Omega,\mathcal F_0,\wP;\,H^{\gamma- 2/p}_{p,\,\theta +2-p}(\mathcal O)\big),
\end{align*}
and for $p\in[2,\infty)$,
\begin{align*}
\mathfrak H^\gamma_{p,\theta}(\mathcal O,T)&:=\bigg\{u\in\HH^\gamma_{p,\theta-p}(\mathcal O,T)\,:\,u(0,\,\cdot\,)\in U^\gamma_{p,\theta}(\mathcal O) \text{ and } du=f\,dt+\sum_{\kappa=1}^\infty g^\kappa\,dw_t^\kappa\\
& \qquad\qquad\qquad\qquad\qquad\quad\text{for some } f\in \HH^{\gamma-2}_{p,\theta+p}(\mathcal O,T),\;g\in \HH^{\gamma-1}_{p,\theta}(\mathcal O,T;\ell_2)\bigg\},
\end{align*}
equipped with the norm
\[\|u\|_{\mathfrak H^\gamma_{p,\theta}(\mathcal O,T)}:=\|u\|_{\HH^\gamma_{p,\theta-p}(\mathcal O,T)}+\|f\|_{\HH^{\gamma-2}_{p,\theta+p}(\mathcal O,T)}+\|g\|_{\HH^{\gamma-1}_{p,\theta}(\mathcal O,T;\ell_2)}+\|u(0,\,\cdot\,)\|_{U^\gamma_{p,\theta}(\mathcal O)}.\]
The equality $du=f\,dt+\sum_{\kappa=1}^\infty g^\kappa\,dw_t^\kappa$ above is  a shorthand  for
\begin{equation}\label{solNotion}
\langle u(t,\,\cdot\,),\varphi\rangle=\langle u(0,\,\cdot\,),\varphi\rangle +\int_0^t\langle f(s,\,\cdot\,),\varphi\rangle\,ds+\sum_{\kappa=1}^\infty \int_0^t\langle g^\kappa(s,\,\cdot\,),\varphi\rangle\,dw^\kappa_s
\end{equation}
for all $\varphi\in C_0^\infty(\mathcal O),\;t\in[0,T]$.
\end{definition}

\begin{rem}
{\bf(a)}
If $p\in[2,\infty)$, then the sum of stochastic integrals in \eqref{solNotion} converges in the space $\mathcal M_T^{2,c}(\R,(\mathcal F_t))$ of continuous, square integrable, $\R$-valued martingales \frenchspacing{w.r.t} $(\mathcal F_t)_{t\in[0,T]}$.  For the convenience of the readers we \fli include \ilf a proof in Appendix \ref{ConSumInt}.

{\bf(b)} Using the arguments of \textsc{Krylov} in \cite[Remark 3.3]{Kry99}, we get the uniqueness (up to indistinguishability) of the pair $(f,g)\in \HH^{\gamma-2}_{p,\theta+p}(\mathcal O,T)\times\HH^{\gamma-1}_{p,\theta}(\mathcal O,T;\ell_2)$ which fulfils \eqref{solNotion}. Consequently, the norm in $\mathfrak H^\gamma_{p,\theta}(\mathcal O,T)$ is well defined.
\end{rem}

\begin{definition}
We call a predictable $\mathcal D'(\domain)$-valued stochastic process $u=(u(t,\,\cdot\,))_{t\in[0,T]}$ a  \emph{solution}  of equation \eqref{eq} if it is a solution of equation \eqref{solNotion} where $f$ is replaced by $\sum_{\mu,\nu=1}^da^{\mu\nu}u_{x_\mu x_\nu}$ and $u(0,\,\cdot\,)=u_0$.
\end{definition}
The next result is taken from \textsc{Kim} \cite{Kim08}.

\begin{thm}\label{kim}
Let $p\in[2,\infty)$ and $\gamma\in\R$. There exists a constant $\kappa_0\in(0,1)$, depending only on $d,\;p,\;(a^{\mu\nu})_{1\leq\mu,\nu\leq d}$ and $\mathcal O$, such that \fli for any $\theta\in(d-\kappa_0,d-2+p+\kappa_0)$, $g\in\HH_{p,\theta}^{\gamma-1}(\mathcal O,T;\ell_2)$ and $u_0\in U^\gamma_{p,\theta}(\mathcal O)$ the \ilf equation \eqref{eq} has a unique solution in the class $\mathfrak H^\gamma_{p,\theta}(\mathcal O,T)$.\\
For this solution
\begin{equation}\label{kimEstimate}
\|u\|_{\mathfrak H^\gamma_{p,\theta}(\mathcal O,T)}^p\leq C\left(\|g\|_{\HH^{\gamma-1}_{p,\theta}(\mathcal O,T;\ell_2)}^p+\|u_0\|_{U^\gamma_{p,\theta}(\mathcal O)}^p
\right),
\end{equation}
where the constant $C$ depends only on $d,\;p,\;\gamma,\;\theta,\;\fli(a^{\mu\nu})_{1\leq\mu,\nu\leq d}\ilf,\;T$ and $\mathcal O$.
\end{thm}
We will need the following straightforward consequence of this  Theorem \ref{kim}. Recall that if $m\in\Nn$ and $f\in\mathcal D'(\domain)$ is sufficient regular, then $|D^mf|_{\ell_p}$ stands for $(\sum_{|\alpha|=m}|D^\alpha f|^p)^{1/p}$, the (pointwise) $\ell_p$-norm of the vector of the $m$-th order derivatives of $f$.
\begin{cor}\label{cor}
In the situation of Theorem \ref{kim} with $\gamma=m\in\mathds N$, the following inequality holds for every $\tau\in[0,p]$.
\[\int_\Omega\int_0^T\|\rho^{m-\delta}|D^m u(\omega,t,\,\cdot\,)|_{\ell_p}\|_{\Lp(\mathcal O)}^\tau\,dt\,\wP(d\omega) \leq C\left(\|g\|_{\HH^{m-1}_{p,\theta}(\mathcal O,T;\ell_2)}+\|u_0\|_{U^m_{p,\theta}(\mathcal O)}
\right)^\tau,\]
where $\delta=1+\frac{d-\theta}p$.
\end{cor}

\begin{proof}
Theorem \ref{kim} implies, in particular, that
\[\|u\|_{\HH_{p,\theta-p}^m(\mathcal O,T)}\leq C\left(\|g\|_{\HH^{m-1}_{p,\theta}(\mathcal O,T)}+\|u_0\|_{U^m_{p,\theta}(\mathcal O)}
\right),\]
and we have
\begin{align*}
\|u\|_{\HH_{p,\theta-p}^m(\mathcal O,T)}^p&=\int_\Omega\int_0^T
\|u(\omega,t,\,\cdot\,)\|_{H^m_{p,\theta-p}(\mathcal O,T)}^p\,dt\,\wP(d\omega)\\
& \geq C\,\int_\Omega\int_0^T
\sum_{k=0}^m\|\rho^{k+(\theta-p-d)/p}|D^ku(\omega,t,\,\cdot\,)|_{\ell_p}\|_{\Lp(\mathcal O)}^p\,dt\,\wP(d\omega)\\
&\geq C\,\int_\Omega\int_0^T
\|\rho^{m-\delta}|D^m u(\omega,t,\,\cdot\,)|_{\ell_p}\|_{\Lp(\mathcal O)}^p\,dt\,\wP(d\omega)
\end{align*}
with $\delta=1+\frac{d-\theta}p\in\big((2-\kappa_0)/p\;,(p+\kappa_0)/p\big)$.
Now let $\tau\in[0,p]$. Jensen's inequality  for concave functions, see, e.g.\ \textsc{Schilling} \cite[Theorem 12.14]{MIMS},  yields
\begin{align*}
\int_\Omega\int_0^T\|\rho^{m-\delta}|D^m u(\omega,t,\,\cdot\,)|_{\ell_p}\|_{\Lp(\mathcal O)}^\tau\,dt\,\wP(d\omega) &\leq  C(T) \left(\|g\|_{\HH^{m-1}_{p,\theta}(\mathcal O,T)}^p+\|u_0\|_{U^m_{p,\theta}(\mathcal O)}^p
\right)^{\tau/p}\\
&\leq C\left(\|g\|_{\HH^{m-1}_{p,\theta}(\mathcal O,T)}+\|u_0\|_{U^m_{p,\theta}(\mathcal O)}
\right)^\tau.
\end{align*}
 In the last step we have used the fact that all norms on $\R^2$ are equivalent.
\end{proof}

\begin{rem}\label{remDaPra1}
Consider the Hilbert space case $p=2$ and assume $g\in\HH^\gamma_{2,\theta}(\domain,T;\ell_2)$. The expression $\sum_{\kappa=1}^\infty\int_0^tg^\kappa(s,\,\cdot\,)\,dw^\kappa_s$ can be considered as an $\HH^\gamma_{2,\theta}(\domain)$-valued stochastic integral $\int_0^tG(s)\,dW_s$ with respect to a cylindrical Wiener process $(W_t)_{t\in[0,T]}$ on $\ell_2$ whose coordinate processes are $(w^\kappa_t)_{t\in[0,T]},\;\kappa\in\Nn$. (See, \frenchspacing{e.g.} \textsc{Da Prato, Zabczyk} \cite{DaPraZab} or \textsc{Peszat, Zabczyk} \cite{PesZab} for stochastic integration w.r.t. cylindrical processes.) Here $(G(t))_{t\in[0,T]}$ is a stochastic process in the space of Hilbert-Schmidt operators $L_{(\text{HS})}(\ell_2,H^\gamma_{2,\theta}(\domain))$ defined by
\[G(\omega,t):\ell_2\to H^\gamma_{2,\theta}(\domain),\;(x^\kappa)_{\kappa\in\Nn}\mapsto \sum_{\kappa\in\Nn} g^\kappa(\omega,t,\,\cdot\,)x^\kappa,\quad (\omega,t)\in\Omega_T,\]
and it is an element of the space $L_2(\Omega_T; L_{(\text{HS})}(\ell_2,H^\gamma_{2,\theta}(\domain)))$. Indeed, for fixed $(\omega,t)\in\Omega_T$ we have
\begin{align*}
\|G(\omega,t)\|^2_{L_{(\text{HS})}(\ell_2,H^\gamma_{2,\theta}(\domain))} &= \sum_{\kappa\in\Nn}\|g^\kappa(\omega,t,\,\cdot\,)\|^2_{H^\gamma_{2,\theta}(\domain)}\\
&=\sum_{\kappa\in\Nn}\sum_{n\in\Z}c^{n\theta}\|\zeta_{-n}(c^n\,\cdot\,)g^\kappa(\omega,t,c^n\,\cdot\,)\|_{H^\gamma_2(\R^d)}^2\\
&=\sum_{n\in\Z}c^{n\theta}\big\|\zeta_{-n}(c^n\,\cdot\,)g(\omega,t,c^n\,\cdot\,)\big\|_
{H^\gamma_2(\Rd;\ell_2)}^2
\end{align*}
by Tonelli's theorem, so that
\[\|G\|_{L_2(\Omega_T;L_{(\text{HS})}(\ell_2;H^\gamma_{2,\theta}(\domain)))}=\|g\|_{\HH^\gamma_{2,\theta}(\domain,T;\ell_2)}.\]
As a consequence, equation \eqref{eq} can be rewritten in the form
\begin{equation}\label{eq'}
du=\sum_{\mu,\nu=1}^d a^{\mu\nu}u_{x_\mu x_\nu}\,dt + dM_t,\quad u(0,\,\cdot\,)=u_0,
\end{equation}
where $(M_t)_{t\in[0,T]}\in\mathcal M^{2,c}_T(H^\gamma_{2,\theta}(\domain),(\mathcal F_t))$ is the $H^\gamma_{2,\theta}(\domain)$-valued, square-integrable martingale given by
\[M_t:=\int_0^tG(s)\,dW_s,\qquad t\in[0,T].\]
\end{rem}

\indent
\begin{rem}\label{remDaPra2}
In Examples \ref{ex1}, \ref{ex2} and \ref{ex3} below the solution $u$ of equation \eqref{eq} in $\mathfrak H^\gamma_{2,\theta}(\domain,T)$ as given by Theorem \ref{kim} coincides with the weak solution of equation \eqref{eq'} with \fli zero \ilf Dirichlet boundary condition in the sense of \textsc{Da Prato, Zabczyk} \cite{DaPraZab}.

In the examples we consider equation \eqref{eq'} driven by certain Wiener processes $(M_t)_{t\in[0,T]}$ in $L_2(\domain)$ with $u_0\in U^2_{2,2}(\domain),\;d=2$ and the solution $u$ is in the class $\mathfrak H^2_{2,2}(\domain,T)\subset\HH^2_{2,0}(\domain,T)$. (Strictly speaking,  in Example \ref{ex3} $(M_t)_{t\in[0,T]}$ is not a Wiener process, but it is one conditioned on the family of random variables $Y_\lambda,\;\lambda\in\nabla$.) Thus, by Remark \ref{Kufner} (a) we know that $u$ is an element of $L_2(\Omega_T;\,\accentset{\circ}{W}^1_2(\domain))$.
Let us now introduce the operator
\[(A,D(A)):=\Bigg(\sum_{\mu,\nu=1}^d a^{\mu\nu}\frac{\partial^2}{\partial x_\mu\partial x_\nu}\;,\;\bigg\{u\in\accentset{\circ}{W}^1_2(\domain)\,:\,\sum_{\mu,\nu=1}^d a^{\mu\nu}u_{x_\mu x_\nu}\in L_2(\domain)\bigg\}\Bigg)\]
and consider the equation
\begin{equation}\label{eqA}
du(t,\,\cdot\,)=Au(t,\,\cdot\,)\;dt+dM_t,\quad u(0,\,\cdot\,)=u_0\in L_2(\domain),\;t\in[0,T].
\end{equation}
A weak solution of equation \eqref{eqA} in the sense of \textsc{Da Prato, Zabczyk} \cite{DaPraZab} is an $L_2(\domain)$-valued predictable process $u=(u(t,\,\cdot\,))_{t\in[0,T]}$ with $\wP$-almost surely Bochner integrable trajectories $t\mapsto u(\omega,t,\,\cdot\,)$ satisfying
\begin{equation}\label{solDaZab}
\langle u(t,\,\cdot\,),\zeta\rangle_{L_2(\domain)}=\langle u_0,\zeta\rangle_{L_2(\domain)}+\int_0^t\langle u(s,\,\cdot\,),A^*\zeta\rangle_{L_2(\domain)}\,ds +\langle M_t,\zeta\rangle_{L_2(\domain)}
\end{equation}
for all $t\in[0,T]$ and $\zeta\in D(A^*)$. It is given by the variation of constants formula
\[u(t,\,\cdot\,)=e^{tA}u_0+\int_0^te^{(t-s)A}\,dM_s,\quad t\in[0,T],\]
where $(e^{tA})_{t\geq0}$ is the contraction semigroup on $L_2(\domain)$ generated by $A$.

It is clear that the solution $u\in\mathfrak H^2_{2,2}(\domain,T)$ given by Theorem \ref{kim} satisfies
\begin{equation}\label{solDaZab'}
\langle u(t,\,\cdot\,),\varphi\rangle_{L_2(\domain)}=\langle u_0,\varphi\rangle_{L_2(\domain)}+\int_0^t\langle A^{1/2}u(s,\,\cdot\,),A^{1/2}\varphi\rangle_{L_2(\domain)}\,ds +\langle M_t,\varphi\rangle_{L_2(\domain)}
\end{equation}
for all $t\in[0,T]$ and $\varphi\in C_0^\infty(\domain)$.  Note that the operator $A$ is self-adjoint because the coefficients $a^{\mu\nu},\;1\leq\mu,\nu\leq d$, are constants.  Since every $\zeta\in D(A^*)=D(A)\subset\accentset{\circ}{W}^1_2(\domain)$ is the limit in $W^1_2(\domain)$ of a sequence of test functions $(\varphi_k)_{k\in\Nn}\subset C_0^\infty(\domain)$, one can go to the limit $k\to\infty$ for $\varphi=\varphi_k$ in \eqref{solDaZab'} to obtain equation \eqref{solDaZab}.
\end{rem}

\section{Besov regularity  for  SPDEs}\label{resultSection}
In this section we state and prove our main result.
\fli We give some concrete examples to illustrate its applicability\ilf.
The  result is formulated in terms of the $L_\tau$-spaces
\[L_{\tau}(\Omega_T;B^s_{\tau,\tau}(\domain))=L_{\tau}(\Omega_T,\mathcal P,\wP\otimes\lambda;B^s_{\tau,\tau}(\domain)),\quad \tau\in (0,\infty),\;s\in(0,\infty),\]
and the spaces introduced in the last section.


\begin{thm}\label{result}
Let $p\in[2,\infty)$ \fli and \ilf $g\in \HH^{\gamma-1}_{p,\theta}(\mathcal O,T;\ell_2),\;u_0\in U^\gamma_ {p,\theta}(\mathcal O)$ for some $\gamma\in\Nn$ and $\theta\in(d-\kappa_{0},d-2+p+\kappa_{0})$  with $\kappa_0=\kappa_0(d,p,(a^{\mu\nu}),\mathcal O)\in(0,1)$ as in Theorem \ref{kim}.
Let $u$ be the unique solution in the class $\mathfrak H^\gamma_{p,\theta}(\mathcal O,T)$ of equation \eqref{eq} and assume furthermore that
\begin{equation}\label{addass}
 u\in \Lp\big(\Omega_T;\,B^s_{p,p}(\domain)\big)\quad\text{ for some }\quad s\in\Big(0,\gamma\wedge\fli\Big(1+\frac{d-\theta}p\Big)\ilf\Big].
\end{equation}
Then, \fli we have
\[u\in L_\tau(\Omega_T;B^\alpha_{\tau,\tau}(\domain)),\quad\frac 1\tau=\frac \alpha d+\frac 1p,\quad\text{for all }\quad\alpha\in\Big(0,\gamma\wedge\frac{sd}{d-1}\Big),\]
and the following estimate holds\ilf
\begin{equation}\label{mainIneq}
\|u\|_{L_\tau(\Omega_T;B^\alpha_{\tau,\tau}(\domain))}\leq C\left(\|g\|_{\HH^{\gamma-1}_{p,\theta}(\mathcal O,T;\ell_2)}+\|u_0\|_{U^\gamma_{p,\theta}(\mathcal O)}
+\|u\|_{\Lp(\Omega_T;B^s_{p,p}(\domain))}\right).
\end{equation}
\fli Here the constant $C$ depends only on $d,\;p,\;\gamma,\;\alpha,\;s,\;\theta,\;(a^{\mu\nu})_{1\leq\mu,\nu\leq d},\;T$ and $\mathcal O$. \ilf
\end{thm}

\begin{rem}
Since the constant $\kappa_0=\kappa_0(d,p,\fli(a^{\mu\nu})\ilf,\mathcal O)$ is greater than zero, we can always choose $\theta=d$. In this case, we know from Theorem \ref{kim} that for each $\gamma\in\Nn$ we have a unique solution $u$ in the class $\mathfrak H^\gamma_{p,d}(\domain,T)$, provided the free term $g$ and the initial condition $u_0$ \fli are sufficiently regular\ilf. In particular, we get
\begin{equation*}
u\in\HH^{\gamma}_{p,d-p}(\domain,T)=L_{p}(\Omega_T,\mathcal P,\wP\otimes\lambda;H^\gamma_{p,d-p}(\domain))\subseteq L_p(\Omega_T;W^1_{p}(\domain))\subseteq L_{p}(\Omega_T;B^1_{p,p}(\domain)).
\end{equation*}
Thus, the additional requirement \eqref{addass} is fulfilled with $s=1$.
Since $\domain$ is  an arbitrary bounded Lipschitz domain \fli it is in general not clear if $u$ belongs to $L_p(\Omega_T;W^s_p(\mathcal O))$ for all $s<2$, \ilf
compare Example \ref{ex1} below.

However, if $\gamma \geq 2$ our result shows that we obtain higher regularity than $s=1$ in \fli the \ilf nonlinear approximation scale, namely
\begin{equation*}
 u\in\fli L_\tau\ilf(\Omega_T;B^\alpha_{\tau,\tau}(\domain)) ,\quad\frac{1}{\tau}=\frac{\alpha}{d}+\frac{1}{p},\quad \text{ for all }\quad\alpha<\frac{d}{d-1}.
\end{equation*}
\end{rem}


\begin{proof}[Proof of Theorem \ref{result}]
We fix $\alpha$ and $\tau$ as stated in the theorem and choose a wavelet Riesz-basis
\begin{equation*}
\left\{\varphi_k,\;\psi_\ijk\,:\,(\ijk)\in\{1,\cdots,2^d-1\}\times\N\times\Z^d\right\}
\end{equation*}
of $L_2(\Rd)$ which \fli fulfils \ilf the \fli assumptions \ilf from Section \ref{waveletsBesovSpaces} with \fli$r>\gamma$\ilf.
Given $(j,k)\in\N\times\Zd$ let
\begin{equation*}
Q_\jk:=2^{-j}k+2^{-j}\,[-N,N]^{d},
\end{equation*}
such that $\fli\supp\psi_\ijk\subset Q_\jk\ilf$ for all $i\in\{1,\ldots,2^d-1\}$ and $\supp\varphi_k\subset Q_{0,k}$ for all $k\in\Zd$. Remember that the supports of the corresponding dual basis \fli fulfil \ilf the same requirements.
For our purpose the set of all indices associated with that wavelets that may have common support with the domain $\domain$ will play an important role and we denote them by
\begin{equation*}
 \Lambda:=\big\{(\ijk)\in\{1,\ldots,2^d-1\}\times\N\times\Zd\,\big|\,Q_{\jk}\cap\domain\neq\emptyset\big\}.
\end{equation*}
In particular, we will also use the following notation:
\begin{equation*}
 \Gamma:=\{k\in\Zd:Q_{0,k}\cap\mathcal O\neq\emptyset\}.
\end{equation*}

Due to the assumption $u\in \Lp(\Omega_T;B^s_{p,p}(\domain))$ we have $u(\omega,t,\,\cdot\,)\in B^s_{p,p}(\domain)$ for $\wP\otimes\lambda$-almost every $(\omega,t)\in\Omega_T$. As $\mathcal O$ is a Lipschitz domain there exists a linear and bounded extension operator $\mathcal E:B^s_{p,p}(\domain)\to B^s_{p,p}(\Rd)$, i.e. there exists a constant $C>0$ such that for $\wP\otimes\lambda$-almost every $(\omega,t)\in\Omega_T$:
\begin{equation*}
\mathcal E u(\omega,t,\,\cdot\,)\big|_{\mathcal O} = u(\omega,t,\,\cdot\,)\qquad \text{and} \qquad
\|\mathcal E u(\omega,t,\,\cdot\,)\|_{B^s_{p,p}(\Rd)} \leq C  \|u(\omega,t,\,\cdot\,)\|_{B^s_{p,p}(\domain)},
\end{equation*}
see, e.g.\ \textsc{Rychkov} \cite{Rych}.
In the \fli sequel \ilf we will  omit  the $\mathcal E$ in our notation and write $u$ instead of $\mathcal E u$.

Theorem \ref{BesovChar01} tells us that for almost all $(\omega,t)\in\Omega_T$ the following equality holds on the domain $\mathcal O$
\begin{align*}
u(\omega,t,\,\cdot\,)&=\sum_{k\in\Gamma}\langle u(\omega,t,\,\cdot\,),\widetilde\varphi_k\rangle\varphi_k + \sum_{(\ijk)\in\Lambda}\langle u(\omega,t,\,\cdot\,),\widetilde\psi_{\ijk,p'}\rangle\psi_{\ijk,p},
\end{align*}
where the sums converge unconditionally in $B^s_{p,p}(\Rd)$. Furthermore,  cf.\  Corollary \ref{BesovChar02}, we get for $\wP\otimes\lambda$-almost all $(\omega,t)\in\Omega_T$
\begin{equation}\label{A}
\|u(\omega,t,\,\cdot\,)\|_{B^\alpha_{\tau,\tau}(\domain)}^\tau
\leq C \Big(\sum_{k\in\Gamma}|\langle u(\omega,t,\,\cdot\,),\widetilde\varphi_k\rangle|^\tau+ \sum_{(\ijk)\in\Lambda}|\langle u(\omega,t,\,\cdot\,),\widetilde\psi_{\ijk,p'}\rangle|^\tau\Big).
\end{equation}
Hence,  it is enough to prove  that
\begin{equation}\label{B}
\int_\Omega\int_0^T \sum_{k\in\Gamma}|\langle u(\omega,t,\,\cdot\,),\widetilde\varphi_k\rangle|^\tau\,dt\,\wP(d\omega)
	\,\leq\, C\, \|u\|_{L_p(\Omega_T; B^s_{p,p}(\domain))}^\tau
\end{equation}
\fli and \ilf
\begin{equation}\label{C}
\begin{aligned}
\lefteqn{\int_\Omega\int_0^T\sum_{(\ijk)\in\Lambda}|\langle u(\omega,t,\,\cdot\,),\widetilde\psi_{i,j,k,p'}\rangle|^\tau\,dt\,\wP(d\omega)}\\
&\qquad\qquad\qquad\leq C\left(\|g\|_{\HH^{\gamma-1}_{p,\theta}(\mathcal O,T;\ell_2)}+\|u_0\|_{U^\gamma_{p,\theta}(\mathcal O)}
+\|u\|_{\Lp(\Omega_T;\,B^s_{p,p}(\domain))}\right)^\tau.
\end{aligned}
\end{equation}

We start with \eqref{B}. The index set $\Gamma$ introduced above is finite because of the boundedness of the domain $\domain$, so that we can use Jensen's inequality to get for $\wP\otimes\lambda$-almost all $(\omega,t)\in\Omega_T$
\begin{align*}
 \sum_{k\in\Gamma} |\langle u(\omega,t,.),\widetilde{\varphi}_k\rangle|^\tau
	\leq C\bigg(\Big( \sum_{k\in\Gamma} |\langle u(\omega,t,.),\widetilde{\varphi}_k\rangle|^p\Big)^{1/p} \bigg)^\tau	
	\leq C\,\|u(\omega,t,.)\|_{B^s_{p,p}(\domain)}^\tau.
\end{align*}
In the last step we used Theorem \ref{BesovChar01} and the boundedness of the extension operator. Integration with respect to $\wP\otimes\lambda$ and another application of Jensen's inequality yield \eqref{B}.

Now let us focus on  the  inequality \eqref{C}. To this end, we introduce the following notation
\begin{align*}
\rho_\jk&:=\text{dist}(Q_{j,k},\partial\domain)=\inf_{x\in Q_\jk}\rho(x),\\
\Lambda_j &:= \big\{ (i,l,k)\in\Lambda\,:\,l=j \big\} ,\\
\Lambda_{j, m}&:=\left\{(\ijk)\in\Lambda_j\,:\,\; m2^{-j}\leq\rho_\jk<( m+1)2^{-j}\right\}	 ,\\
\Lambda_j^0&:=\Lambda_j\setminus\Lambda_{j,0},\\
\Lambda^0&:=\bigcup_{j\in\N}\Lambda_j^0,
\end{align*}
where $j, m\in\N$ and $k\in\Zd$.
We split the expression on the left hand side of $\eqref{C}$ into
\begin{align}\label{split}
\int_\Omega\int_0^T \sum_{(\ijk)\in\Lambda^0}
	& \big| \langle u(\omega,t,\,\cdot\,),\widetilde\psi_{i,j,k,p'} \rangle \big|^\tau\,dt\,\wP(d\omega) \nonumber\\
	& +\int_\Omega\int_0^T\sum_{(\ijk)\in\Lambda\setminus\Lambda^0} \big| \langle u(\omega,t,\,\cdot\,),\widetilde\psi_{i,j,k,p'} \rangle \big|^\tau \,dt\,\wP(d\omega)
	=: I + II
\end{align}
and estimate each term separately.

Let us begin with $I$. Fix $(\ijk)\in\Lambda^0$ and $(\omega,t)\in\Omega_T$ such that
\[\int_\mathcal O\big|\rho(x)^{\gamma-s}|D^\gamma u(\omega,t,x)|_{\ell_p}\big|^p\,dx<\infty.\]
 By  Corollary \ref{cor} this \fli holds \ilf for $\wP\otimes\lambda$-almost all $(\omega,t)\in\Omega_T$. By a Whitney-type inequality, \fli also \ilf known as the Deny-Lions lemma, see, \frenchspacing{e.g.} \textsc{DeVore}, \textsc{Sharpley} \cite[Theorem 3.4]{DeVSh84},  there exists a polynomial $P_{j,k}$ of total degree less than $\gamma$ such that
\[\|u(\omega,t,\,\cdot\,)-P_{j,k}\|_{\Lp(Q_\jk)}\leq C2^{-j\gamma}|u(\omega,t,\,\cdot\,)|_{W^\gamma_p(Q_\jk)},\]
where the last norm is finite since $\rho_{j,k}=\text{dist}(Q_\jk,\partial \mathcal O)>0$. Since $\widetilde\psi_{\ijk,p'}$ is orthogonal to every polynomial of total degree less than $\gamma$, one gets
\begin{align*}
\big|\langle u(\omega,t,\,\cdot\,),\widetilde\psi_{\ijk,p'}\rangle\big|
	&=\big|\langle u(\omega,t,\,\cdot\,)-P_{\jk},\widetilde\psi_{\ijk,p'}\rangle\big|\\
	&\leq\|u(\omega,t,\,\cdot\,)-P_{\jk}\|_{\Lp(Q_\jk)} \, \|\widetilde\psi_{\ijk,p'}\|_{L_{p'}(Q_\jk)}\\
	&\leq C \, 2^{-j\gamma}\,\big|u(\omega,t,\,\cdot\,)\big|_{W^\gamma_p(Q_\jk)}\\
	&\leq C \, 2^{-j\gamma}\rho_\jk^{s-\gamma}\Big(\int_{Q_\jk}\big|\rho(x)^{\gamma-s}\,|D^\gamma u(\omega,t,x)|_{\ell_p}\big|^p\,dx\Big)^{1/p}\\
	&=:C \, 2^{-j\gamma}\rho_\jk^{s-\gamma} \, \mu_\jk(\omega,t).
\end{align*}
 Fix  $j\in\N$. Summing over all indices $(\ijk)\in\Lambda_j^0$ and applying H{\"o}lder's inequality with exponents  $\frac p\tau>1$ and $\frac p{p-\tau}$ one finds
\begin{align}\label{Hoelder}
\sum_{(\ijk)\in\Lambda_j^0} \big| \langle u(\omega,t,\,\cdot\,),\widetilde\psi_{\ijk,p'} \rangle \big|^\tau
	&\leq C \sum_{(\ijk)\in\Lambda_j^0} 2^{-j\gamma\tau}\rho_\jk^{(s-\gamma)\tau}\mu_\jk(\omega,t)^\tau	\nonumber\\
	&\leq C \Big( \sum_{(\ijk)\in\Lambda_j^0}\mu_\jk(\omega,t)^p \Big)^{\frac\tau p}
		\Big( \sum_{(\ijk)\in\Lambda_j^0}2^{\frac{-pj\gamma\tau}{p-\tau}} \rho_{j,k}^{\frac{(s-\gamma)p\tau}{p-\tau}}\Big)^{\frac{p-\tau}{p}}.
\end{align}
Since \fli any \ilf $x\in \mathcal O$ lies outside of all but at most a constant number $C>0$ of the cubes $Q_\jk$, $k\in\Zd$, we get the following bound for the first factor on the right hand side
\begin{align}\label{Hoelder01}
\bigg( \sum_{(\ijk)\in\Lambda_j^0} \mu_\jk(\omega,t)^p \bigg)^{\frac\tau p}
	&=\bigg( \sum_{(\ijk)\in\Lambda_j^0}\int_{Q_\jk}\big|\rho(x)^{\gamma-s}\,|D^\gamma u(\omega,t,x)|_{\ell_p}\big|^p\,dx\bigg)^{\frac\tau p} \nonumber\\
	&\leq C \,\left\|\rho^{\gamma-s}|D^\gamma u(\omega,t,\,\cdot\,)|_{\ell_p}\right\|_{\Lp(\mathcal O)}^\tau.
\end{align}
In order to estimate the second factor \fli in \eqref{Hoelder} \ilf we use the Lipschitz character of the domain $\domain$ which implies that
\begin{align}\label{LipFin}
 |\Lambda_{j, m}|\leq C 2^{j(d-1)}\qquad\text{ for all $j, m\in\N$.}
\end{align}
 The constant $C>0$ does not depend on $j$ or $m$. Moreover, the  boundedness  of $\domain$ yields  $\Lambda_{j, m}=\emptyset$ for all $j, m\in\N$ with $ m\geq C2^j$. \fli Consequently,
\begin{equation}\label{uglySums}
\begin{aligned}
\bigg( \sum_{(\ijk)\in\Lambda_j^0} 2^{\frac{-pj\gamma\tau}{p-\tau}} \rho_{j,k}^{\frac{(s-\gamma)p\tau}{p-\tau}}\bigg)^{\frac{p-\tau}p}
	&\leq \bigg( \sum_{ m=1}^{C 2^j} \sum_{(\ijk)\in\Lambda_{j, m}} 2^{\frac{-pj\gamma\tau}{p-\tau}} \rho_{j,k}^{\frac{(s-\gamma)p\tau}{p-\tau}}\bigg)^{\frac{p-\tau}p}	\\
	&\leq C \bigg( \sum_{ m=1}^{C2^j} 2^{j(d-1)} \, 2^{-j\frac{p\gamma\tau}{p-\tau}} ( m \, 2^{-j})^{\frac{(s-\gamma)p\tau}{p-\tau}}\bigg)^{\frac{p-\tau}p}	\\
	&\leq C \bigg( 2^{j\left(d-1-\frac{s p \tau}{p-\tau}\right)}+ 2^{j\left(d-\frac{\gamma p \tau}{p-\tau}\right)} \bigg)^{\frac{p-\tau}p}.
\end{aligned}
\end{equation}
Now, let us sum over all $j\in\N$ and integrate over $\Omega_T$ with respect to $\wP\otimes\lambda$ on both sides of \fli the \ilf inequality \eqref{Hoelder}. By using \eqref{uglySums} and \eqref{Hoelder01} from above and Corollary \ref{cor} we get
\begin{align*}
 \int_\Omega&\int_0^T\sum_{(\ijk)\in\Lambda^0} \big| \langle u(\omega,t,\,\cdot\,),\widetilde\psi_{i,j,k,p'} \rangle \big|^\tau\,dt\,\wP(d\omega)	\\
	&\leq 	C
	 	\sum_{j\in\N} \bigg( 2^{j\left(d-1-\frac{s p \tau}{p-\tau}\right)} + 2^{j\left(d-\frac{\gamma p \tau}{p-\tau}\right)} \bigg)^\frac{p-\tau}{p}
		\int_\Omega\int_0^T \left\| \rho^{\gamma-s}\,|D^\gamma u(\omega,t,\,\cdot\,)|_{\ell_p} \right\|_{\Lp(\domain)}^\tau \,dt\,\wP(d\omega)\\
	&\leq 	C
		\bigg(
			\sum_{j\in\N}  2^{j\left(d-1-\frac{s p \tau}{p-\tau}\right)\frac{p-\tau}{p} }
		      + \sum_{j\in\N} 2^{j\left(d-\frac{\gamma p \tau}{p-\tau}\right)\frac{p-\tau}{p} }
		\bigg)
		\left(\|g\|_{\HH^{\gamma-1}_{p,\theta}(\mathcal O,T;\ell_2)}+\|u_0\|_{U^\gamma_{p,\theta}(\mathcal O)}\right)^\tau.
\end{align*}
One can see that the sums on the right hand side converge  if, and only if, $\alpha\in\left(0,\gamma\wedge s\frac{d}{d-1}\right)$. \ilf Finally,
\begin{align}\label{D}
 \int_\Omega\int_0^T \sum_{(\ijk)\in\Lambda^0} \big| \langle u(\omega,t,\,\cdot\,),\widetilde\psi_{i,j,k,p'} \rangle \big|^\tau\,dt\,\wP(d\omega)
	\leq C\left(\|g\|_{\HH^{\gamma-1}_{p,\theta}(\mathcal O,T;\ell_2)}+\|u_0\|_{U^\gamma_{p,\theta}(\mathcal O)}\right)^\tau.
\end{align}

%

Now we estimate the term $II$ in \eqref{split}. \fli First we fix $j\in\N$ and use H{\"o}lder's inequality and
\eqref{LipFin} to get
\begin{align*}
 \sum_{(\ijk)\in\Lambda_{j,0}} \big| \langle u(\omega,t,.),\widetilde{\psi}_{\ijk,p'} \rangle \big|^\tau
	&\leq C\,2^{j(d-1)\frac{p-\tau}{p}}
		\Big( \sum_{(\ijk)\in\Lambda_{j,0}} \big| \langle u(\omega,t,.),\widetilde{\psi}_{\ijk,p'} \rangle \big|^p \Big)^\frac{\tau}{p}.
\end{align*}
Summing over all $j\in\N$ and using H{\"o}lder's inequality again, yields
\begin{align*}
 \sum_{(\ijk)\in\Lambda\backslash\Lambda^0}& \big| \langle u(\omega,t,.),\widetilde{\psi}_{\ijk,p'} \rangle \big|^\tau
	=	\sum_{j\in\N} \Big[\sum_{(\ijk)\in\Lambda_{j,0}} \big| \langle u(\omega,t,.),\widetilde{\psi}_{\ijk,p'} \rangle \big|^\tau\Big] \\
	&\leq 	C \sum_{j\in\N} \Big[ 2^{j(d-1)\frac{p-\tau}{p}} \Big( \sum_{(\ijk)\in\Lambda_{j,0}} \big| \langle u(\omega,t,.),\widetilde{\psi}_{\ijk,p'} \rangle \big|^p \Big)^\frac{\tau}{p} \Big] \\
	&\leq	C \bigg(\sum_{j\in\N} 2^{j\left( \frac{(d-1)(p-\tau)}{p}-sp \right)\frac{p}{p-\tau}}\bigg)^{\frac{p-\tau}{p}}
		\bigg( \sum_{j\in\N} \sum_{(\ijk)\in\Lambda_{j,0}} 2^{j s p}\big| \langle u(\omega,t,.),\widetilde{\psi}_{\ijk,p'} \rangle \big|^p \bigg)^\frac{\tau}{p}.
\end{align*}\ilf
Using Theorem \ref{BesovChar01} and the  boundedness  of the extension operator, one gets for $\wP\otimes\lambda$-almost every $(\omega,t)\in\Omega_T$ that
\begin{align*}
 \sum_{(\ijk)\in\Lambda\backslash\Lambda^0}& \big| \langle u(\omega,t,.),\widetilde{\psi}_{\ijk,p'} \rangle \big|^\tau
	\leq	C \| u(\omega,t,.) \|_{B^s_{p,p}(\domain)}^{\tau}
		\bigg(\sum_{j\in\N} 2^{j\left( \frac{(d-1)(p-\tau)}{p}-sp \right)\frac{p}{p-\tau}}\bigg)^{\frac{p-\tau}{p}}.
\end{align*}
The series on the right hand side converges if and only if $\alpha\in\left( 0,s\frac{d}{d-1} \right)$. But this is \fli part of our assumptions\ilf, so that for $\wP\otimes\lambda$-almost every $(\omega,t)\in\Omega_T$
\begin{align*}
 \sum_{(\ijk)\in\Lambda\backslash\Lambda^0}& \big| \langle u(\omega,t,.),\widetilde{\psi}_{\ijk,p'} \rangle \big|^\tau
	\leq C \| u(\omega,t,.) \|_{B^s_{p,p}(\domain)}^{\tau}.
\end{align*}
Let us integrate over $\Omega_T$ with respect to $\wP\otimes\lambda$ and use Jensen's inequality to get \fli
\begin{align*}
\int_\Omega\int_0^T\!\!
	\sum_{(\ijk)\in\Lambda\setminus\Lambda^0} \big| \langle u(\omega,t,\,\cdot\,), \widetilde\psi_{i,j,k,p'} \rangle \big|^\tau dt \,\wP(d\omega)
	&\leq C\int_\Omega\int_0^T\|u(\omega,t,\,\cdot\,)\|_{B^s_{p,p}(\domain)}^\tau\,dt\,\wP(d\omega)\\
	&\leq C\left[\int_\Omega\int_0^T\|u(\omega,t,\,\cdot\,)\|_{B^s_{p,p}(\domain)}^p\,dt\,\wP(d\omega)\right]^{\frac\tau p}.
\end{align*}
Because of \eqref{D} this proves \eqref{C}. Now \eqref{B} and \eqref{A} finish the proof.\ilf
\end{proof}


Next, we give some examples for an application of Theorem \ref{result}. \fli We \ilf are mainly interested in the Hilbert space case $p=2$.

\begin{ex}\label{ex1}
Let us first consider equation \eqref{eq} in the form \eqref{eq'} where the driving process $(M_t)_{t\in[0,T]}$ is a Wiener process in $\fli\accentset{\circ}{W}^1_2(\mathcal O)\ilf$ with covariance operator $Q\in L_{(\text{nuc})}(\fli\accentset{\circ}{W}^1_2(\mathcal O)\ilf)$. It can be represented as  a  stochastic integral process $(\int_0^tG(s)\,dW_s)_{t\in[0,T]}$ \frenchspacing{w.r.t.} the cylindrical Wiener process $(W_t)_{t\in[0,T]}$ on $\ell_2$ by defining the integrand process $(G(t))_{t\in[0,T]}$ in the space of Hilbert-Schmidt operators $L_{(\text{HS})}(\ell_2,\fli\accentset{\circ}{W}^1_2(\mathcal O)\ilf)$ as the constant deterministic process
\begin{equation}\label{G}
G(\omega,t):\ell_2\to \fli\accentset{\circ}{W}^1_2(\domain)\ilf,\;(x^\kappa)_{\kappa\in\Nn}\mapsto\fli\sum\nolimits_{\kappa\in\Nn}\ilf\sqrt{\lambda_\kappa}x^\kappa e_\kappa,\quad (\omega,t)\in\Omega_T,
\end{equation}
where $(e_\kappa)_{\kappa\in\Nn}$ is an orthonormal basis of $\fli\accentset{\circ}{W}^1_2(\domain)\ilf$ consisting of eigenvectors of $Q$ \fli with positive eigenvalues $(\lambda_\kappa)_{\kappa\in\Nn}$\ilf.

This corresponds to defining $g=(g^\kappa)_{\kappa\in\Nn}$ in equation \eqref{eq} by
\begin{equation}\label{g}
g^\kappa(\omega,t,\,\cdot\,):=\sqrt{\lambda_\kappa}e_\kappa,\qquad \kappa\in\Nn,\;(\omega,t)\in\Omega_T.
\end{equation}
It is easy to see that $g$ is an element of $\HH^1_{2,d}(\domain,T;\ell_2)$. By definition
\begin{equation}
\begin{aligned}\label{estg}
\|g\|_{\HH^1_{2,d}(\domain,T;\ell_2)}^2&= T^2\sum_{n\in\Z}c^{nd}\big\|\zeta_{-n}(c^n\,\cdot\,)(\sqrt{\lambda_\kappa}e_\kappa(c^n\,\cdot\,))_{\kappa\in\Nn}\big\|^2_{H^1_2(\Rd;\ell_2)}\\
&=T^2\sum_{\kappa\in\Nn}\lambda_\kappa\sum_{n\in\Z}c^{nd}\|\zeta_{-n}(c^n\,\cdot\,)e_\kappa(c^n\,\cdot\,)\|_{H^1_2(\Rd)}^2\\
&=T^2\sum_{\kappa\in\Nn}\lambda_\kappa\|e_\kappa\|_{H^1_{2,d}(\domain)}^2.
\end{aligned}
\end{equation}
 Using the norm equivalence  \eqref{ne}, one \fli has
\begin{align*}
\|g\|_{\HH^1_{2,d}(\domain,T;\ell_2)}^2\,\leq\, CT^2\sum_{\kappa\in\Nn}\lambda_\kappa\sum_{|\alpha|\leq 1}\|\rho^{|\alpha|}D^\alpha e_\kappa\|_{L_2(\domain)}^2
&\,\leq\, CT^2\sum_{\kappa\in\Nn}\lambda_\kappa\sum_{|\alpha|\leq 1}\|D^\alpha e_\kappa\|_{L_2(\domain)}^2\\
&\,=\,C T^2\sum_{\kappa\in\Nn}\lambda_\kappa\;<\;\infty.
\end{align*}\ilf
Thus, in a 2-dimensional setting, Theorem \ref{kim} with $d=\theta=\gamma=2$ tells us that for every initial condition $u_0\in U^2_{2,2}(\mathcal O)= L_2(\Omega,\mathcal F_0,\wP;H^1_{2,2}(\domain))$ equation $\eqref{eq}$ has a unique solution $u$ in the class $\mathfrak H^2_{2,2}(\domain,T)\subset\HH^2_{2,0}(\domain,T)=L_2(\Omega_T;H^2_{2,0}(\domain))$. As a trivial consequence,
\[u\in L_2(\Omega_T;W^1_2(\domain))=L_2(\Omega_T;B^1_{2,2}(\domain))\]
because we have the equality
\[H^2_{2,0}(\domain)=\big\{u\in\mathcal D'(\domain)\,:\,\rho^{|\alpha|-1}D^\alpha u \in L_2(\domain)
\text{ for all }\alpha\in\N^2\text{ with }|\alpha|\leq2\big\}.\]
(In fact, according to Remark \ref{Kufner} we even know that $u\in L_2(\Omega_T;\accentset\circ W^1_2(\domain))$.)

Note that in general $u$ does not belong to $L_2(\Omega_T;W^s_2(\mathcal O))$ \fli for all $s<2$.
Since $\mathcal O$ is an arbitrary bounded Lipschitz domain, certain second derivatives might explode near the boundary \ilf and the norm $\|u(\omega,t,\,\cdot\,)\|_{W^2_2(\domain)}$ as well as $\|u(\omega,t,\,\cdot\,)\|_{W^s_2(\domain)}$, where  $s\in(1,2)$, might not be finite. If $\mathcal O$ is a polygonal domain, one can derive an explicit upper bound for the regularity in the Sobolev scale $L_2(\Omega_T;W^s_2(\mathcal O)),\;s>0$. \fli Adapting \ilf techniques used in \fli \textsc{Grisvard} \cite{Gri85}, \cite{Gri92} \ilf to our stochastic setting, one can show that $u\notin L_2(\Omega_T;\,W^s_2(\domain))$ if $s>1+\pi/\gamma_0$, where $\gamma_0$ is the measure of the largest interior angle at a corner of $\partial O$.


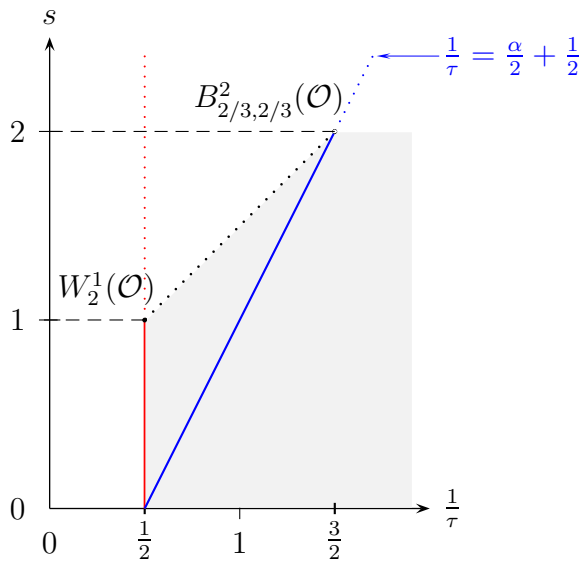
\begin{figure}

\psset{yunit=2.5,xunit=2.5}

\begin{center}
\begin{pspicture}(0,0)(2,3)	
	\psaxes[linewidth=0.5pt,arrowscale=1.5]{->}(2,2.5)
	\uput[r](2,0){$\frac{1}{\tau}$}
	\uput[u](0,2.5){$s$}
	
	\psline(0.5,-0.03)(0.5,0)
	\psline[linecolor=red](0.5,0)(0.5,1)
	\psline[linecolor=red,linestyle=dotted](0.5,1)(0.5,2.4)
	
	 \pspolygon[linecolor=mylightGray,fillcolor=mylightGray,fillstyle=solid](0.51,0.01)(0.51,0.99)(1.5,1.99)(1.9,1.99)(1.9,0.01)
	
	\psline[linecolor=blue](0.5,0)(1.5,2)
	\psline[linecolor=blue,linestyle=dotted](1.5,2)(1.7,2.4)
	
	\psline[linewidth=1pt,linestyle=dotted](0.5,1)(1.5,2)

	\psdots[dotscale=.5](.5,1)
	\uput[d](0.5,0){$\frac{1}{2}$}
	\psline[linewidth=0.4pt,linestyle=dashed](-0.03,1)(0.5,1)
	\uput{.2}[u](0.3,1){$W^{1}_2(\domain)$}
	
	
	\psdots[dotscale=.5,dotstyle=o](1.5,2)
	\psline[linewidth=.4pt,linestyle=dashed](0,2)(1.5,2)
	\uput[d](1.5,0){$\frac{3}{2}$}
	\psline(1.5,-0.03)(1.5,0.03)
	\uput[ul](1.6,2){$B_{2/3,2/3}^{2}(\domain)$}

	\uput[r](2,2.4){\rnode[linecolor=blue]{G}{\blue $\frac{1}{\tau}=\frac{\alpha}{2}+\frac{1}{2}$}}
	\pnode(1.7,2.4){H}
	 \nccurve[nodesepA=.05,nodesepB=.05,linecolor=blue,linewidth=0.3pt,angleA=180,angleB=0]{->,arrowlength=2.5}{G}{H}
	
\end{pspicture}

\end{center}
\caption[DeVoreTriebel002]{\begin{tabular}[t]{l} \pci Besov regularity in the scale $B^\alpha_{\tau,\tau}(\domain)$, \fli$1/\tau=\alpha/2+1/2,$ \ilf vs. Sobolev regularity \icp \\ \pci of the solution, illustrated in a \textsc{DeVore-Triebel} diagram. \icp  \end{tabular}}\label{fig:DeVoreTriebel002}

\end{figure}

However, in the considered situation Theorem \ref{result} with $s=1$ states that we have
\[u\in L_\tau(\Omega_T;B^\alpha_{\tau,\tau}(\domain))\]
for every $\alpha<2$ and $1/\tau=\alpha/2+1/2$. This constellation is illustrated in Figure \ref{fig:DeVoreTriebel002}, where each point $(1/ \tau, s)$ represents the smoothness spaces of functions with ``$ s$ derivatives in $L_{\tau}(\domain)$''. Based on the knowledge \fli that \ilf $u\in L_2(\Omega_T;\,W^1_2(\domain))$ and $u\in L_\tau(\Omega_T;\,B^\alpha_{\tau,\tau}(\domain))$ for all $\alpha<2,\;1/\tau=\alpha/2+1/2$,  interpolation and embedding theorems yield that $u$ also belongs to each of the spaces \fli$L_{\tau}(\Omega_T;\,B^s_{\tau,\tau}(\domain)),\;0<\tau<2,\;s<(1/2+1/\tau)\wedge2$\ilf. This is indicated by the shaded area.
\end{ex}

\begin{ex}\label{ex2}
In view of equality \eqref{estg} it is clear that we can apply Theorem \ref{kim} and Theorem \ref{result} in the same \fli way \ilf as in Example \ref{ex1}, \frenchspacing{i.e.} with $d=\theta=\gamma=2$ and $s=1$, \fli if the driving process $(M_t)_{t\in[0,T]}$ in \eqref{eq'} is a Wiener process in $W^1_2(\domain)$ with covariance operator $Q\in L_{(\text{nuc})}(W^1_2(\domain))$, and even if it is a Wiener process in $H^1_{2,2}(\domain)$ with covariance operator $Q\in L_{(\text{nuc})}(H^1_{2,2}(\domain))$. In the first case $(M_t)_{t\in[0,T]}$ does not satisfy a zero Dirichlet boundary condition as in Example \ref{ex1}, and in the second case $(M_t)_{t\in[0,T]}$ behaves even more irregularly near the boundary in the sense that the first derivatives are allowed to blow up near $\partial\domain$.

In these cases we choose $(e_\kappa)_{\kappa\in\Nn}$ in \eqref{G} and \eqref{g} to be an orthonormal basis of the space $W^1_2(\domain)$, respectively $H^1_{2,2}(\domain)$, consisting of eigenvectors of $Q\in L_{(\text{nuc})}(W^1_2(\domain))$, respectively $Q\in L_{(\text{nuc})}(H^1_{2,2}(\domain))$, with corresponding eigenvalues $(\lambda_\kappa)_{\kappa\in\Nn}$.\ilf

As in Example \ref{ex1} the solution $u$ lies in $L_2(\Omega_T;W^1_2(\domain))=L_2(\Omega_T;B^1_{2,2}(\domain))$ and, by Theorem \ref{result}, it also lies in $L_2(\Omega_T;B^\alpha_{\tau,\tau}(\domain)),\;1/\tau=\alpha/2+1/2,\;\alpha<2$, see Figure \ref{fig:DeVoreTriebel002}.
\end{ex}

\begin{ex}\label{ex3}
\fli Let the driving process $(M_t)_{t\in[0,T]}$ in \eqref{eq'} be \ilf a time-dependent version of the stochastic wavelet  expansion  introduced in \textsc{Abramovich} et \frenchspacing{al.} \cite{Abr} in the context of Bayesian nonparametric regression and generalized in \textsc{Bochkina} \cite{Boch}, \textsc{Cioica} et \frenchspacing{al.} \cite{Rit}. This noise model is formulated in terms of a wavelet basis expansion on the domain $\domain\subset\Rd$ with random coefficients of \fli prescribed \ilf sparsity and thus \fli tailor-made \ilf for applying adaptive techniques with regard to the numerical approximation of the corresponding SPDEs. Via the choice of certain \fli parameters \ilf specifying the distributions of the wavelet coefficients it also allows for an explicit control of the spatial Besov regularity of $(M_t)_{t\in[0,T]}$ . We first describe the general noise model and then deduce a further example for the application of Theorem \ref{result}.

Let $\{\psi_\lambda\,:\,\lambda\in\nabla\}$ be a multiscale Riesz basis for $L_2(\domain)$ consisting of scaling functions at a fixed scale level $j_0\in\Z$ and of wavelets at level $j_0$ and all finer levels. As in the introduction, the notation we use here is different from that used in Section \ref{waveletsBesovSpaces} because we do not consider a basis on the whole space $\Rd$ but on the bounded domain $\domain$. Information like scale level, spatial location and type of the wavelets or scaling functions are encoded in the indices $\lambda\in\nabla$.  We refer to \textsc{Cohen} \cite[Sections 2.12, 2.13 and 3.9]{C03} and \textsc{Dahmen, Schneider} \fli \cite{DS98}, \cite{DS99a}, \cite{DS99b} \ilf for detailed descriptions of multiscale bases on bounded domains. Adopting the notation of \textsc{Cohen} we write $\nabla=\bigcup_{j\geq j_0-1}\nabla_j$, where for $j\geq j_0$ the set $\nabla_j\subset\nabla$ contains the indices of all wavelets $\psi_\lambda$ at scale level $j$ and where $\nabla_{j_0-1}\subset\nabla$ is the index set referring to the scaling functions at scale level $j_0$ which we denote by $\psi_\lambda,\;\lambda\in\nabla_{j_0-1}$, for the sake of notational simplicity. We make the following assumptions concerning our basis. Firstly, the cardinalities of the index sets $\nabla_j,\;j\geq j_0-1$, satisfy
\begin{equation}\label{lim}
 C^{-1}2^{jd}\leq  |\nabla_j| \leq C  2^{jd},\qquad j\geq j_0-1.
\end{equation}
Secondly, we assume that the basis admits norm equivalences similar to those described in Theorem \ref{BesovChar01}. There exists an $r\in\Nn$ (depending on the smoothness of the scaling functions $\psi_\lambda,\;\lambda\in\nabla_{j_0-1}$, and on the degree of polynomial exactness of their linear span), such that, given $p,q>0,\;\max\{0,d(1/p-1)\}<s<r$, and a real valued distribution $f\in\mathcal D'(\domain)$, we have $f\in B^s_{p,q}(\domain)$ if and only if $f$ can be represented as $f=\sum_{\lambda\in\nabla}c_{\lambda}\psi_{\lambda},\;(c_\lambda)_{\lambda\in\nabla}\subset\R$ (convergence in $\mathcal D'(\domain)$), \fli such that
\begin{equation}\label{equi}
\Bigg(\sum_{j={j_0-1}}^\infty 2^{jq(s+d(\frac 12-\frac 1p))}
\Big(\sum_{\lambda\in\nabla_j}|c_\lambda|^p\Big)^{q/p}\Bigg)^{1/q}<\infty.
\end{equation}
Furthermore, $\|f\|_{B^s_{p,q}(\domain)}$ is equivalent to the quasi-norm \eqref{equi}. \ilf
Concrete constructions of bases satisfying these assumptions can be found in the  literature mentioned above.
Concerning the family of independent standard Brownian motions $(w^\kappa_t)_{t\in[0,T]},\;\kappa\in\Nn$, in \eqref{eq} respectively \eqref{eq'}, we modify our notation and write $(w_t^\lambda)_{t\in[0,T]},\;\lambda\in\nabla,$ instead. The  description  of the noise model involves \fli parameters \ilf $a\geq0,\;b\in[0,1],\;c\in\R$, with $a+b>1$. For every $j\geq j_0-1$ we set $\sigma_j=(j-(j_0-2))^{\frac{c d}2}2^{-\frac{a (j-(j_0-1))d}2}$ and let $Y_\lambda,\;\lambda\in\nabla_j,$ be Bernoulli distributed random variables on $(\Omega,\mathcal F_0,\wP)$ with parameter $p_j=2^{-b(j-(j_0-1))d}$, such that the random variables and processes $Y_\lambda,\;(w_t^\lambda)_{t\in[0,T]},\;\lambda\in\nabla,$ are stochastically independent. Now we are ready to define $(M_t)_{t\in[0,T]}$ by
\begin{equation}\label{model}�
M_t:=\sum_{j=j_0-1}^\infty\sum_{\lambda\in\nabla_j}
\sigma_jY_{\lambda}\psi_{\lambda}\cdot w^\lambda_t,\qquad t\in[0,T].
\end{equation}
Using \eqref{equi}, \eqref{lim} and $a+b>1$, it is easy to check that the infinite sum converges in $L_2(\Omega_T;L_2(\mathcal O))$ as well as in the space $\mathcal M_T^{2,c}(L_2(\domain),(\mathcal F_t))$ of continuous, square integrable, $L_2(\domain)$-valued martingales \frenchspacing{w.r.t.} the filtration $(\mathcal F_t)_{t\in[0,T]}$. Moreover, by the choice of the hyperparameters $a,\;b$ and $c$ one has an explicit control of the convergence of the infinite sum in \eqref{model} in the (quasi-)Banach spaces $L_{p_2}(\Omega_T;B^s_{p_1,q}(\domain)),\;s<r,\;p_1,q>0,\;p_2\leq q$. (Compare \textsc{Cioica} et \frenchspacing{al.} \fli \cite{Rit} which \ilf can easily be adapted to our setting.)

With regard to Theorems \ref{kim} and \ref{result} let again $d=p=\gamma=\theta=2$. Equation \eqref{eq'} with $(M_t)_{t\in[0,T]}$ defined as above corresponds to equation \eqref{eq} if we set \[g^{\lambda}(\omega,t,\,\cdot\,):=\sigma_jY_\lambda(\omega)\psi_\lambda(\,\cdot\,),\qquad \lambda\in\nabla_j,\;j\geq j_0-1,\;(\omega,t)\in\Omega_T,\]
and sum over all $\lambda\in\nabla$ instead of $\kappa\in\Nn$. In the following we write $\ell_2=\ell_2(\nabla)$. \fli Since $a+b>1$ and $\|g\|_{\HH_{2,2}^0(\domain,T;\ell_2)}=\sqrt{2/T}\|M\|_{L_2(\Omega_T;L_2(\mathcal O))}$ we have $g\in\HH_{2,2}^0(\domain,T;\ell_2)$\ilf. Let us impose a bit more smoothness on $g$ and assume that $a+b>2$. This is sufficient to ensure that $g\in\HH^1_{2,2}(\mathcal O,T;\ell_2)$: Using \eqref{ne} one sees that the $\HH^1_{2,2}(\mathcal O,T;\ell_2)$-norm of $g=(g^\lambda)_{\lambda\in\nabla}$ satisfies
\begin{align*}
\|g\|_{\HH^1_{2,2}(\mathcal O,T;\ell_2)}^2
&=\E\int_0^T\sum_{n\in\Z}c^{n2}\|\zeta_{-n}(c^n\,\cdot\,)g(t,c^n\,\cdot\,)\|^2_{H^1_2(\Rd;\ell_2)}\,dt\\
&=\E\int_0^T\sum_{\lambda\in\nabla}\|g^\lambda(t,\,\cdot\,)\|_{H^1_{2,2}(\domain)}^2\,dt\\
&=T\,\E\sum_{j=j_0-1}^\infty\sum_{\lambda\in\nabla_j}\sigma_j^2Y_\lambda^2\|\psi_\lambda\|_{H^1_{2,2}(\domain)}^2\\
&\leq C \sum_{j=j_0-1}^\infty\sum_{\lambda\in\nabla_j}\sigma_j^2p_j\sum_{|\alpha|\leq 1}\|\rho^{|\alpha|}D^\alpha \psi_\lambda\|_{L_2(\domain)}^2\\
&\leq C \sum_{j=j_0-1}^\infty\sum_{\lambda\in\nabla_j}\sigma_j^2p_j\|\psi_\lambda\|_{W^1_2(\domain)}^2.
\end{align*}
Since $W^1_2(\domain)=B^1_{2,2}(\domain)$ with equivalent norms we can use the equivalence \eqref{equi} with $f=\psi_\lambda$ to get
\begin{align*}
\|g\|_{\HH^1_{2,2}(\mathcal O,T;\ell_2)}^2
&\leq C \sum_{j=j_0-1}^\infty\sum_{\lambda\in\nabla_j}\sigma_j^2p_j2^{2j}\\
&=C\sum_{j=j_0-1}^\infty|\nabla_j|(j-(j_0-2))^{2c}2^{-2 a(j-(j_0-1))}2^{-2 b(j-(j_0-1))}2^{2j}\\
&\leq C \sum_{j=j_0-1}^\infty (j-(j_0-2))^{2c}2^{-2j(a+b-2)}.
\end{align*}
In  the last step  we used  \eqref{lim}  with $d=2$. Thus $g\in\HH^1_{2,2}(\mathcal O,T;\ell_2)$. As in Example \ref{ex1} we may apply Theorems \ref{kim} and \ref{result} to conclude that for every initial condition $u_0\in L_2(\Omega,\mathcal F_0,\wP; H^1_{2,2}(\domain))$ there exists a unique solution of equation \eqref{eq} in the class $\mathfrak H^2_{2,2}(\domain,T)$, which, in general, is not in $L_2(\Omega_T,\fli W^s_2(\domain))\ilf$ for all $\fli s<2\ilf$, but it belongs to every space $L_2(\Omega_T;B^\alpha_{\tau,\tau}(\domain))$ with $\alpha<2$ and \fli$\tau=2/(\alpha+1)$\ilf.
\end{ex}

\begin{rem}
\pci
In practice, many adaptive wavelet-based algorithms are  realized with the energy norm of the problem which is equivalent to a Sobolev norm. Let us denote by $\{\eta_\lambda\,:\,\lambda\in\nabla\}$ a wavelet Riesz basis of $W^s_2(\domain)$ for some $s>0$, which can be obtained by rescaling the wavelet basis $\{ \psi_\lambda\,:\,\lambda\in\nabla\}$ of $L_2(\domain)$, see, e.g.\ \textsc{Cohen} \cite{C03} or \textsc{Dahmen} \cite{D97}.
For the best $N$-term approximation in this Sobolev norm, it is well known that
\begin{equation*}
 u\in B^{\alpha}_{\tau,\tau}(\domain),\;\frac{1}{\tau}=\frac{\alpha-s}{d}+\frac{1}{2}\quad \Longrightarrow\quad \sigma_{N,W^s_2(\domain)}(u) \leq C\,N^{-(\alpha-s)/d},
\end{equation*}
where
\begin{equation*}
\sigma_{N,W^s_2(\domain)}(u):=\inf\Big\{\nnrm{u-u_N}{W^s_2(\domain)}\,:\, u_N=\sum_{\lambda\in\Lambda}c_\lambda\eta_\lambda\;:\;\Lambda\subset\nabla,\;|\Lambda|\leq N,\;c_\lambda\in\R,\;\lambda\in\Lambda\Big\}.
\end{equation*}
Therefore, similar to the $L_2(\domain)$-setting, the approximation order of \fli the \ilf best $N$-term wavelet scheme in $W^s_2(\domain)$ depends on the Besov regularity of the object one wants to approximate.

There exist adaptive wavelet-based algorithms which are guaranteed to converge and which indeed asymptotically realize the convergence rate of best $N$-term approximation with respect to the Sobolev norm. For example, \textsc{Cohen}, \textsc{Dahmen}, \textsc{DeVore} \cite{CDD1} designed such an adaptive numerical scheme for solving (deterministic) elliptic PDEs. First results for parabolic problems were obtained by \textsc{Schwab}, \textsc{Stevenson} \cite{SchSt09}.

Once again, the use of adaptive algorithms is justified if the rate of approximation that can be achieved is higher \fli than in \ilf classical uniform schemes. Let $u_N$, $N\in\Nn$, denote a uniform approximation scheme (e.g.\ a Galerkin approximation) of $u$. It is well-known that under certain natural conditions, see, e.g.\ \textsc{Dahlke}, \textsc{Dahmen}, \textsc{DeVore} \cite{DDD} or \textsc{DeVore} \cite{DeV98} or \textsc{Hackbusch} \cite{Hack96},
\begin{equation*}
 \nnrm{u-u_N}{W^{s}_2(\domain)}\leq C N^{-(\alpha-s)/d}\nnrm{u}{W^{\alpha}_{2}(\domain)}.
\end{equation*}
\fli This means that, even \ilf in this case, adaptivity can pay off if the Besov smoothness of the solution is higher than its Sobolev regularity.

Let us discuss this relationship in more detail for the examples from above. We consider approximation in $W^1_2(\domain)$.  As already mentioned in Example \ref{ex1}, in general we \fli cannot \ilf expect that the spatial Sobolev regularity of the solution is higher than $3/2$. Therefore, uniform schemes yield an approximation rate of $O(N^{-1/4})$.


\begin{figure}

\psset{yunit=2.5,xunit=2.5}

\begin{center}
\begin{pspicture}(0,0)(2,3)	
	\psaxes[linewidth=0.5pt,arrowscale=1.5]{->}(2,2.5)
	\uput[r](2,0){$\frac{1}{\tau}$}
	\uput[u](0,2.5){$s$}
	
	\pspolygon[linecolor=mylightGray,fillcolor=mylightGray,fillstyle=solid](0.51,0.01)(0.51,1.49)(1.5,1.99)(1.9,1.99)(1.9,0.01)

	\psline(0.5,-0.03)(0.5,0)
	\psline[linecolor=red](0.5,0)(0.5,1.5)
	\psline[linecolor=red,linestyle=dotted](0.5,1.5)(0.5,2.4)
	
	\psline[linecolor=blue](0.5,0)(1.5,2)
	\psline[linecolor=blue,linestyle=dotted](1.5,2)(1.6,2.2)

 	\psline[linecolor=ddblau](0.5,1)(0.833,1.666)
	\psline[linecolor=ddblau,linestyle=dotted](0.833,1.666)(1.2,2.4)
	
	\psline[linecolor=magenta,linestyle=dotted](0.5,1.5)(1.5,2)
	

	\psdots[dotscale=.5](.5,1.5)
	\uput[d](0.5,0){$\frac{1}{2}$}
	\uput[l](-0.03,1.5){$\frac{3}{2}$}
	\psline[linewidth=0.4pt](-0.03,1.5)(0.03,1.5)
	\uput{.4}[l](-0.1,1.25){\rnode{C}{$W^{3/2}_2(\domain)$}}
	\pnode(0.5,1.5){D}
	\nccurve[linewidth=.3pt,angleA=0,angleB=240,nodesep=.1]{->}{C}{D}

	\psdots[dotscale=.5,dotstyle=o](1.5,2)
	\uput[d](1.5,0){$\frac{3}{2}$}
	\psline(1.5,-0.03)(1.5,0.03)
	\uput{.4}[r](1.6,2){\rnode{E}{$B_{2/3,2/3}^{2}(\domain)$}}
	\pnode(1.5,2){F}
	\nccurve[linewidth=.3pt,angleA=180,angleB=0,nodesep=.1]{->}{E}{F}

	\uput{1}[r](1.2,2.4){\rnode[linecolor=ddblau]{G}{{\color{ddblau} $\frac{1}{\tau}=\frac{\alpha-1}{2}+\frac{1}{2}$}}}
	\pnode(1.2,2.4){H}
	\nccurve[nodesepA=.05,nodesepB=.05,linecolor=ddblau,linewidth=0.3pt,angleA=180,angleB=0]{->,arrowlength=2.5}{G}{H}
	
	\uput{1}[r](1,1){\rnode[linecolor=blue]{G}{\blue $\frac{1}{\tau}=\frac{\alpha}{2}+\frac{1}{2}$}}
	\pnode(1,1){H}
	\nccurve[nodesepA=.05,nodesepB=.05,linecolor=blue,linewidth=0.3pt,angleA=180,angleB=0]{->,arrowlength=2.5}{G}{H}

	\psdots[dotscale=.5,dotstyle=o](0.833,1.666)
	\uput{1.8}[r](0.833,1.666){\rnode[linecolor=blue]{I}{$B^{5/3}_{6/5,6/5}(\domain)$}}
	\pnode(0.833,1.666){J}
	 \nccurve[nodesepA=.05,nodesepB=.05,linewidth=0.3pt,angleA=180,angleB=0]{->,arrowlength=2.5}{I}{J}
	
\end{pspicture}

\end{center}
\caption[DeVoreTriebel003]{\begin{tabular}[t]{l} \pci Besov regularity in the scale $B^\alpha_{\tau,\tau}(\domain)$, \fli$1/\tau=(\alpha-1)/2+1/2,$ \ilf  vs. Sobolev \\ regularity of the solution illustrated in a \textsc{DeVore-Triebel} diagram. \icp\end{tabular}}\label{fig:DeVoreTriebel003}

\end{figure}
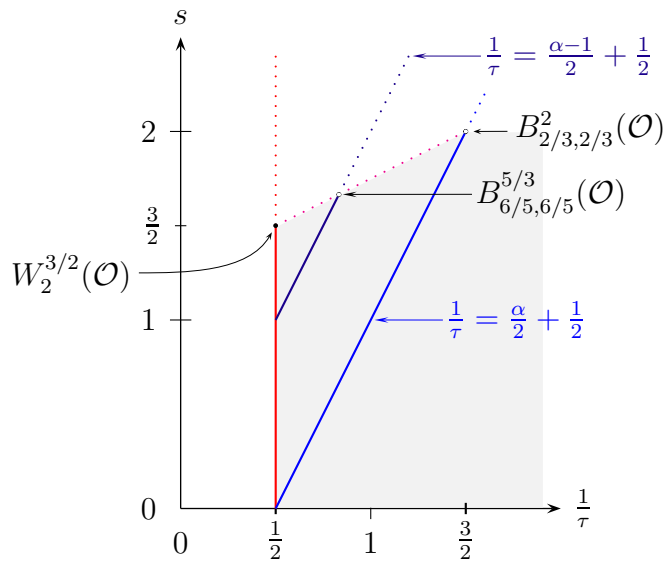

On the other hand, our main result \fli shows that\ilf
\begin{equation*}
 u\in L_{\tau}(\Omega_T;B^\alpha_{\tau,\tau}(\domain)),\;\frac{1}{\tau}=\frac{\alpha}{2}+\frac{1}{2}\;\text{ for all } \alpha<2.
\end{equation*}
Therefore, by interpolation and embedding of Besov spaces we can achieve that the solution is contained in all the spaces $L_\tau(\Omega_T;B^\alpha_{\tau,\tau}(\domain))$ corresponding to the points in the \fli trapezoid with vertices $(1/2,0)$, $(1/2,3/2)$, $(3/2,2)$, $(3/2,0)$ and to the points to the right of this trapezoid \ilf in the \textsc{DeVore-Triebel} diagram, cf. Figure \ref{fig:DeVoreTriebel003}.
As a consequence, we get by a short computation, that
\begin{equation*}
 u\in L_{\tau}(\Omega_T;B^\alpha_{\tau,\tau}(\domain)), \; \frac{1}{\tau}=\frac{\alpha-1}{2}+\frac{1}{2}\;\text{ for all } \alpha<\frac{5}{3}.
\end{equation*}
Thus, best $N$-term wavelet approximation provides order $O(N^{-1/3})$, so that again the use of adaptivity is completely justified.
\icp
\end{rem}

\begin{appendix}
\section{Convergence of the stochastic integrals}\label{ConSumInt}

In this section we give a proof of the $\mathcal M^{2,c}_T(\R,(\mathcal F_t))$-convergence of the sum of the stochastic integral processes $\big(\int_0^t\langle g^\kappa(s,\,\cdot\,),\varphi\rangle\,dw^\kappa_s\big)_{t\in[0,T]},\;\kappa\in\mathds N,$ appearing in formula \eqref{solNotion}. Let us assume that $g\in\HH^\gamma_{p,\theta}(\domain,T;\ell_2)$ for some $p\in[2,\infty)$ and $\gamma,\,\theta\in\R$. We use an  analogous  strategy to \cite{Kry99}, Remark 3.2. Due to the independence of the Brownian motions $(w_t^\kappa)_{t\in[0,T]},\;\kappa\in\mathds N,$ the covariation process $([w^\kappa,w^\ell]_t)_{t\in[0,T]}$ vanishes if $\kappa\neq\ell$, and we have  by It{\^o}'s isometry:
\begin{align*}
\E\left|\sum_{\kappa=1}^\infty \int_0^T\langle g^\kappa(s,\,\cdot\,),\varphi\rangle\,dw^\kappa_s\right|^2
&=\E\left[\sum_{\kappa=1}^\infty \int_0^{(\,\cdot\,)}\langle g^\kappa(s,\,\cdot\,),\varphi\rangle\,dw^\kappa_s\,,\,\sum_{\kappa=1}^\infty \int_0^{(\,\cdot\,)}\langle g^\kappa(s,\,\cdot\,),\varphi\rangle\,dw^\kappa_s\right]_T\\
&=\E\sum_{\kappa=1}^\infty \int_0^T|\langle g^\kappa(s,\,\cdot\,),\varphi\rangle|^2\,ds.
\end{align*}
We are going to show that the last term is less or equal a constant times $\|g\|_{\HH^\gamma_{p,\theta}
(\domain,T;\ell_2)}^2$, which is finite due to our assumption. Then the convergence of the integral processes in $\mathcal M^{2,c}_T(\R,(\mathcal F_t))$ follows by Doob's  maximal inequality  for martingales.

For $u\in \mathcal D'(\mathcal O)$ and $n\in\Z$ we use the notation $u_n:=\zeta_{-n}(c^n\,\cdot\,)u(c^n\,\cdot\,)\in\mathcal S'(\Rd)$. \fli Let us abbreviate $L_\tau(\Rd)$ by $L_\tau$ for all $\tau \geq 1$ in the sequel. Setting $p'=p/(p-1)$, we denote by $\langle\,\cdot\,,\,\cdot\,\rangle_{L_p\times L_{p'}}:L_p\times L_{p'}\to\R$ \ilf the dual form obtained by continuous extension of $\langle\varphi,\psi\rangle=\int \varphi(x)\psi(x)\,dx,\;\varphi,\;\psi\in C_0^\infty(\Rd)$. Now  we are ready to estimate as follows: \fli
\begin{align*}
\lefteqn{\sum_{\kappa=1}^\infty \int\limits_0^T|\langle g^\kappa(s,\,\cdot\,),\varphi\rangle|^2\,ds
\,=\,\sum_{\kappa=1}^\infty \int\limits_0^T\left|\sum_{n\in\Z}c^{nd}\langle g^\kappa_n(s,\,\cdot\,),\varphi_n\rangle\right|^2\,ds}\\
&=\sum_{\kappa=1}^\infty \int\limits_0^T\left|\sum_{n\in\Z}c^{nd}\big\langle(1-\Delta)^{\frac\gamma2} g^\kappa_n(s,\,\cdot\,),(1-\Delta)^{-\frac\gamma2}\varphi_n\big\rangle_{L_p\times L_{p'}}\right|^2\,ds\\
&\leq \sum_{\kappa=1}^\infty \int\limits_0^T\left[\sum_{n\in\Z}c^{nd}\big\||(1-\Delta)^{\frac\gamma2} g^\kappa_n(s,\,\cdot\,)|\cdot|(1-\Delta)^{-\frac\gamma2}\varphi_n|^{1/2}\big\|_{L_2}\cdot
\big\||(1-\Delta)^{-\frac\gamma2}\varphi_n|^{1/2}\big\|_{L_2}\right]^2ds\\
&\leq \sum_{\kappa=1}^\infty \int\limits_0^T\left[\left(\sum_{n\in\Z}c^{2nd}\big\||(1-\Delta)^{\frac\gamma2} g^\kappa_n(s,\,\cdot\,)|\cdot|(1-\Delta)^{-\frac\gamma2}\varphi_n|^{1/2}\big\|_{L_2}^2\right)^{1/2}\right.\\
&\qquad\qquad\qquad\qquad\qquad\qquad\qquad\qquad\qquad\qquad\left.\cdot
\left(\sum_{n\in\Z }\big\|(1-\Delta)^{-\frac\gamma2}\varphi_n\big\|_{L_1}\right)^{1/2}\right]^2\,ds.
\end{align*}\ilf
Here we have used H{\"o}lder's inequality twice.  Since $\varphi$ has compact support in $\domain$ and $\zeta_{-n}$ equals zero outside $\domain_{-n}$,  the functions $\varphi_n$ vanish on $\Rd$ for all but finitely many $n\in\Z$. As a consequence, the sum $\sum_{n\in\Z}\|(1-\Delta)^{-\gamma/2}\varphi_n\|_{L_1}$ has only finitely many non-zero terms. Therefore, \fli
\begin{align*}
\lefteqn{\sum_{\kappa=1}^\infty \int_0^T|\langle g^\kappa(s,\,\cdot\,),\varphi\rangle|^2\,ds}\\
&\leq C\sum_{\kappa=1}^\infty \int_0^T\sum_{n\in\Z}c^{2nd}\big\||(1-\Delta)^{\gamma/2} g^\kappa_n(s,\,\cdot\,)|\cdot|(1-\Delta)^{-\gamma/2}\varphi_n|^{1/2}\big\|_{L_2}^2\,ds\\
&= C\int_0^T\sum_{n\in\Z}c^{2nd}\Big\langle\sum_{\kappa=1}^\infty |(1-\Delta)^{\gamma/2} g^\kappa_n(s,\,\cdot\,)|^2,|(1-\Delta)^{-\gamma/2}\varphi_n|\Big\rangle_{L_{p/2}\times L_{p/(p-2)}}\,ds ,
\end{align*} \ilf
 where the constant $C$ depends on $\varphi$.
In the last step we used the fact that \[(1-\Delta)^{\gamma/2}g_n(s,\,\cdot\,)=\big((1-\Delta)^{\gamma/2}g_n^\kappa(s,\,\cdot\,)\big)_{\kappa\in\Nn}\in L_p(\Rd;\ell_2)\]
 $\wP\otimes\lambda$-almost everywhere in $\Omega_T$, which results from $g$ being an element of $\HH^\gamma_{p,\theta}(\mathcal O;\ell_2)$. Applying again H{\"o}lder's inequality we obtain
\begin{align*}
\lefteqn{\sum_{n\in\Z}c^{2nd}\Big\langle\sum_{\kappa=1}^\infty |(1-\Delta)^{\gamma/2} g^\kappa_n(s,\,\cdot\,)|^2,|(1-\Delta)^{-\gamma/2}\varphi_n|\Big\rangle_{L_{p/2}\times L_{p/(p-2)}}}\\
&\leq \sum_{n\in\Z}c^{2n\theta}\big\||(1-\Delta)^{\gamma/2}g_n(s,\,\cdot\,)|_{\ell_2}\big\|_{L_p}^{1/2}
c^{2n(d-\theta)}\|(1-\Delta)^{-\gamma/2}\varphi_n\|_{L_{p/(p-2)}}\\
&\leq C\left(\sum_{n\in\Z}c^{n\theta}\big\||(1-\Delta)^{\gamma/2}g_n(s,\,\cdot\,)|_{\ell_2}\big\|_{L_p}^p\right)^{\frac 2p}
\left(\sum_{n\in\Z}c^{2n(d-\theta)p/(p-2)}\|(1-\Delta)^{-\gamma/2}\varphi_n\|_{L_{p/(p-2)}}^{p/(p-2)}\right)^{\frac{p-2}p}\\
&\leq C\left(\sum_{n\in\Z}c^{n\theta}\big\||(1-\Delta)^{\gamma/2}g_n(s,\,\cdot\,)|_{\ell_2}\big\|_{L_p}^p\right)^{\frac 2p},
\end{align*}
where we have used that $p\leq 2$ and that only finitely many of the $\varphi_n,\;n\in\Z,$ are non-zero.
All in all we have shown
\[\sum_{\kappa=1}^\infty \int_0^T|\langle g^\kappa(s,\,\cdot\,),\varphi\rangle|^2\,ds\leq C\int_0^T\left(\sum_{n\in\Z}c^{n\theta}\big\||(1-\Delta)^{\gamma/2}g_n(s,\,\cdot\,)|_{\ell_2}\big\|_{L_p}^p\right)^{\frac 2p}ds.\]
Finally, taking the expectation and applying Jensen's inequality
yields
\[\E\sum_{\kappa=1}^\infty \int_0^T|\langle g^\kappa(s,\,\cdot\,),\varphi\rangle|^2\,ds\leq C\|g\|_{\HH^\gamma_{p,\theta}
(\domain,T;\ell_2)}^2\]
and this finishes the proof.


\section{General linear equations}\label{AppB}

In the introduction we have indicated that our main result can be extended to equation \eqref{eq''}. The  major  reason is, that by a result from \textsc{Kim} \cite{Kim08} an estimate similar to the one proved in Corollary \ref{cor} holds not only for the model equation \eqref{eq} but for equations \fli of the type \eqref{eq''},
provided \ilf the coefficients $a^{\mu\nu}, b^\mu$, $c$, $\sigma^{\mu\kappa}$ and $\eta^\kappa$, the free terms $f$ and $g^\kappa$ and the initial value $u_0$ fulfil certain conditions. We can use this fact to extend our regularity result to such equations. In this section we want to get more precise and point out how to do this.

For the  convenience  of the reader we begin by presenting the result from \textsc{Kim} \cite[Theorem 2.8]{Kim08}. Therefore, we need some additional notations.
For $x,y \in \domain$ we shall write $\rho(x,y):=\rho(x)\land\rho(y)$. For $\alpha\in\R$, $\delta \in\left(0,1\right]$ and $k\in\N$ we set:
\begin{align*}
	 [f]_{k}^{(\alpha)} := \sup_{x\in\domain} \rho^{k+\alpha}(x) |D^{k} f(x)| \qquad&\text{ and }\qquad [f]_{k+\delta}^{(\alpha)} := \sup_{\substack{x,y\in\domain \\ |\beta|=k}} \rho^{k+\alpha}(x,y) \frac{|D^{\beta} f(x) - D^{\beta} f(y)|}{|x-y|^{\delta}}, \\
 	|f|_{k}^{(\alpha)} := \sum_{j=0}^{k} [f]_{j}^{(\alpha)} \qquad&\text{ and } \qquad |f|_{k+\delta}^{(\alpha)} := |f|_{k}^{(\alpha)} + [f]_{k+\delta}^{(\alpha)},
\end{align*}
whenever it makes sense. We shall use the same notations for $\ell_2$-valued functions (just replace the absolute values in the above definitions by the $\ell_2$-norms). Furthermore, let's fix an arbitrary function
\begin{equation*}
 \mu_0 : \left[0,\infty\right)\to \left[0,\infty\right),
\end{equation*}
vanishing only on the set of nonnegative integers, i.e. $\mu_0(j)=0$ if and only if $j\in\N$. We set
\begin{equation*}
 t_+:=t+\mu_0(t).
\end{equation*}

Now we are  able  to present the assumptions on the coefficients of equation \eqref{eq''} (see \textsc{Kim} \cite[Assumptions 2.5 and 2.6]{Kim08}).


\begin{itemize}
\item [{\bf [K1]}] For any fixed $x\in\domain$, the coefficients
	\begin{equation*}
	 a^{\mu\nu}\left(.,.,x\right), b^\mu\left(.,.,x\right), c\left(.,.,x\right), \sigma^{\mu\kappa}\left(.,.,x\right), \eta^{\kappa}\left(.,.,x\right) :\Omega\times\left[0,T\right]\to\R
	\end{equation*}
are predictable processes with respect to the given normal filtration $\left(\mathcal F_t\right)_{t\in[0,T]}$.

\item [{\bf [K2]}] (Stochastic parabolicity) There are constants $\delta_0, K > 0$, such that for all $(\omega,t,x)\in\Omega\times\left[0,T\right]\times\domain$ and $\lambda\in\Rd$:
	\begin{equation*}
	 \delta_0 |\lambda|^2 \leq \overline{a^{\mu\nu}}(\omega,t,x)\lambda_\mu\lambda_\nu \leq K|\lambda|^2,
	\end{equation*}
where $\overline{a^{\mu\nu}}:=a^{\mu\nu}-\frac12\langle\sigma^\mu,\sigma^\nu\rangle_{\ell_2}$ for $\mu,\nu\in\{1,\dots,d\}$.

\item [{\bf [K3]}] For all $(\omega,t)\in\Omega\times\left[0,T\right]$:
	\begin{align*}
 	| a^{\mu\nu}(\omega,t,.) |_{|\gamma|_+}^{(0)} + | b^{\mu}(\omega,t,.) |_{|\gamma|_+}^{(1)} &+ | c(\omega,t,.) |_{|\gamma|_+}^{(2)} \\
	 & + | \sigma^{\mu}(\omega,t,.) |_{|\gamma|_+}^{(0)} + | \nu(\omega,t,.) |_{|\gamma+1|_+}^{(1)} \leq K.
	\end{align*}

\item [{\bf [K4]}] The coefficients $a^{\mu\nu}$ and $\sigma^{\mu}$ are uniformly continuous in $x\in\domain$, i.e. for any $\epsilon > 0$ there is a $\delta=\delta(\epsilon)>0$, such that
\begin{equation*}
 | a^{\mu\nu}(\omega,t,x)- a^{\mu\nu}(\omega,t,y)| + | \sigma^\mu(\omega,t,x)-\sigma^\mu(\omega,t,y) |_{\ell_2} \leq \epsilon,
\end{equation*}
for all $(\omega,t)\in\Omega\times\left[0,T\right]$, \fli whenever \ilf $x,y\in\domain$ with $|x-y| \leq \delta$.

\item [{\bf [K5]}] The behaviour of the coefficients $b^\mu$, $c$ and $\nu$ can be controlled near the boundary of $\domain$ in the following way:
\begin{equation*}
 \lim_{\substack{ \rho(x)\to0 \\ x\in\domain }} \sup_{\substack{ \omega\in\Omega \\ t\in[0,T]}}
	\{ \rho(x) | b^\mu(\omega,t,x) | +\rho^2(x) | c(\omega,t,x) | + | \nu(\omega,t,x) |_{\ell_2} \}=0.
\end{equation*}
\end{itemize}

Here is the main result of \textsc{Kim} \cite{Kim08}.

\begin{thm}\label{KimGenerell}
Let $p\in\left[2,\infty\right)$ and let assumptions \emph{[K1]} -- \emph{[K5]} be satisfied with $K,\delta_0 > 0$.
Then there exists a constant
$\kappa_0=\kappa_0(d,p,\delta_0,K,\domain) \in (0,1)$
such that, if
$
 \theta\in(d-\kappa_0,d+\kappa_0+p-2),
$
for any
$
f	\in \HH^{\gamma-2}_{p, \theta+p}(\domain,T)
$,
$
g	\in \HH^{\gamma-1}_{p, \theta}(\domain,T;\ell_2)
$ and
$
	u_0	\in U_{p,\theta}^{\gamma}(\domain),
$
equation \eqref{eq''} with initial value $u_0$ admits a unique solution $u\in \mathfrak{H}^\gamma_{p,\theta}(\domain,T)$,
i.e., there exists an (up to indistinguishability) unique $\mathcal{D}'(\domain)$-valued predictable process
$
 	u\in \HH^{\gamma}_{p, \theta-p}(\domain,T),
$
such that for any $\varphi \in C_0^{\infty}(\domain)$ the equality
\begin{align*}
 \langle u(t,.), \varphi \rangle
	&= \langle u(0,.), \varphi \rangle \\
	&\quad+ \int_{0}^{t} \langle a^{\mu\nu}(s,.)u_{x_\mu x_\nu}(s,.)+b^\mu(s,.) u_{x_\mu}(s,.)+c(s,.)u(s,.)+f(s,.), \varphi \rangle \mathrm{d}s \\
	&\quad+ \sum_{\kappa=1}^{\infty}\int_{0}^{t} \langle \sigma^{\mu\kappa}(s,.)u_{x_\mu}(s,.)+\eta^\kappa(s,.) u(s,.)+g^\kappa(s,.), \varphi \rangle \mathrm{d}w_s^\kappa
\end{align*}
holds for all $t\in\left[0,T\right]$ with probability $1$.
Moreover, for this solutions we have
\begin{equation}\label{kimEstimateGenerell}
 \nnrm{u}{\mathfrak{H}^\gamma_{p,\theta}(\domain,T)}^{p}
	\leq C \left(
		\nnrm{f}{\HH^{\gamma-2}_{p, \theta+p}(\domain,T)}^{p}
		+ \nnrm{g}{\HH^{\gamma-1}_{p, \theta}(\domain,T;\ell_2)}^{p}
		+ \nnrm{u_0}{U_{p,\theta}^{\gamma}(\domain)}^p
		\right),
\end{equation}
where $C$ is a constant depending only on $d$, $\gamma$, $p$, $\theta$, $\delta_0$, $K$, $T$ and the domain $\domain$.
\end{thm}


An immediate consequence of this  theorem  is the following estimate.

\begin{cor}\label{corGenerell}
In the situation of Theorem \ref{KimGenerell} with $\gamma=m\in\mathds N$, the following inequality holds for every $\tau\in[0,p]$:
\begin{align*}
\int_\Omega \int_0^T \|\rho^{m-\beta}|D^m u(\omega,t,\,\cdot\,)|_{\ell_p}\|&_{\Lp(\mathcal O)}^\tau\,dt\,\wP(d\omega)  \nonumber\\
&\leq
C \left(\|f\|_{\HH^{m-2}_{p, \theta+p}}+\|g\|_{\HH^{m-1}_{p,\theta}(\mathcal O,T;\ell_2)}+\|u_0\|_{U^m_{p,\theta}(\mathcal O)}
\right)^\tau,
\end{align*}
where $\beta=1+\frac{d-\theta}p$.
\end{cor}

\begin{proof}
Just repeat the arguments of the proof of Corollary \ref{cor} and use estimate \eqref{kimEstimateGenerell} instead of \eqref{kimEstimate} at the beginning.
\end{proof}


Now we can present our main result in the generalized setting.

\begin{thm}\label{resultGenerell}
Let $p\in[2,\infty)$ and let assumptions \emph{[K1] \textendash{} [K5]} be satisfied with appropriate constants $K,\delta_0 > 0$. Moreover, let $f\in\HH_{p,\theta+p}^{\gamma-2}(\domain,T)$, $g\in \HH^{\gamma-1}_{p,\theta}(\mathcal O,T;\ell_2)$ and $u_0\in U^\gamma_ {p,\theta}(\mathcal O)$ for some $\gamma\in\Nn$. Denote by $u$ the unique solution of equation \eqref{eq''} in the class $\mathfrak{H}_{p,\theta}^{\gamma}(\domain,T)$ for a given $\theta\in(d-\kappa_0,d-2+p+\kappa_0)$, where $\kappa_0=\kappa_0(d,p,\delta_0,K,\domain) \in (0,1)$ as in Theorem \ref{KimGenerell}.
Assume furthermore that $u \in \Lp\big(\Omega_T;\,B^s_{p,p}(\domain)\big)$ for some $s\in(0,\gamma\wedge(1+\frac{d-\theta}p)]$.
\\
Then, \fli we have
\begin{equation*}
 u\in L_\tau(\Omega_T;B^\alpha_{\tau,\tau}(\domain)),\quad\frac 1\tau=\frac \alpha d+\frac 1p,\quad\text{ for all }\quad \alpha\in\Big(0,\gamma\wedge \frac{sd}{d-1}\Big),
\end{equation*}
and the following estimation holds
\begin{align*}\label{mainIneqGenerell}
\|u\|&_{L_\tau(\Omega_T;B^\alpha_{\tau,\tau}(\domain))}
&\leq C \left(\|f\|_{\HH^{\gamma-2}_{p,\theta+p}(\domain,T)} + \|g\|_{\HH^{\gamma-1}_{p,\theta}(\domain,T;\ell_2)} + \|u_0\|_{U^\gamma_{p,\theta}(\domain)} + \|u\|_{\Lp(\Omega_T;B^s_{p,p}(\domain))} \right).
\end{align*}\ilf
\end{thm}

\begin{proof}
We can argue like we did in the proof of Theorem \ref{result}. We just have to use Corollary \ref{corGenerell} where we used Corollary \ref{cor}.
\end{proof}

\end{appendix} 


\section*
\fli
\noindent Petru~A.~Cioica, Stephan~Dahlke, Stefan~Kinzel:
Philipps-Universit{\"a}t Marburg,
FB Mathematik und Informatik, AG Numerik/Optimierung,
Hans-Meerwein-Stra{\ss}e,
35032 Marburg, Germany.\\
\{cioica, dahlke, kinzel\}@mathematik.uni-marburg.de\\

\noindent Felix~Lindner, Ren{\'e}~L.~Schilling:
TU Dresden,
FB Mathematik, Institut f{\"u}r Mathematische Stochastik,
Zellescher Weg 12-14,
01069 Dresden, Germany.\;
\{felix.lindner, rene.schilling\}@tu-dresden.de \\

\noindent Thorsten~Raasch:
Johannes-Gutenberg-Universit{\"a}t Mainz,
Institut f{\"u}r Mathematik, AG Numerische Mathe-matik,
Staudingerweg 9,
55099 Mainz, Germany.\;
raasch@uni-mainz.de\\

\noindent Klaus~Ritter:
TU Kaiserslautern,
Department of Mathematics, Computational Stochastics Group,
Erwin-Schr{\"o}dinger-Stra{\ss}e,
67663 Kaiserslautern, Germany.\;
ritter@mathematik.uni-kl.de
\ilf

\end{document}